\documentclass[a4paper,twoside,10pt]{article}

\usepackage{graphicx}
\usepackage[dvips]{color}
\usepackage[all]{xy}
\usepackage[english]{babel}

\usepackage{amsmath,amssymb}

\usepackage{amsthm}

\usepackage{amscd}

\usepackage{mathrsfs}

\begin{document}

\title{On  the complexity group of stable curves}
\author{Simone Busonero, Margarida Melo, Lidia Stoppino}

\newcommand{\av}{``}
\newcommand{\Z}{\mathbb Z}
\newcommand{\kapp}{\mbox{\itshape {\textbf k}}}
\newcommand{\kappy}{\mbox{\small \itshape {\textbf k}}}
\newcommand{\forn}{\forall n\in\mathbb{N}}
\newcommand{\raw}{\to}
\newcommand{\Raw}{\Rightarrow}
\newcommand{\law}{\gets}
\newcommand{\Law}{\Leftarrow}
\newcommand{\Lra}{\Leftrightarrow}
\newcommand{\spc}{\,\mbox{spec}\,}

\newtheorem{teo}{Theorem}[section]
\newtheorem{prop}[teo]{Proposition}
\newtheorem{lem}[teo]{Lemma}
\newtheorem{rem}[teo]{Remark}
\newtheorem{defi}[teo]{Definition}
\newtheorem{ex}[teo]{Example}
\newtheorem{cor}[teo]{Corollary}
\newtheorem{conj}[teo]{Conjecture}
\pagestyle{myheadings}
\markboth{\small{S. Busonero, M. Melo, L. Stoppino}}{\small{\textit{The complexity of stable curves}}}

\renewcommand{\theequation}{\arabic{section}.\arabic{equation}}

\maketitle
\begin{abstract}
In this paper, we study combinatorial properties of stable curves. 
To the dual graph of any nodal curve, it is naturally associated a group, which is the group of components of the N\'eron model
of the generalized Jacobian of the curve. 
We study the order of this group, called the complexity.
In particular, we provide a partial characterization of the stable curves having maximal complexity, and we 
provide an upper bound, depending only on the genus $g$ of the curve, on the maximal complexity of stable curves; 
this bound is asymptotically sharp for $g\gg 0$. 
Eventually, we state some conjectures on the behavior of stable curves with maximal complexity, and prove partial results in this direction.
\end{abstract}


\section*{Introduction}


Combinatorics is often of great importance in the study of the moduli space of stable curves of genus $g$, $\overline M_g$.
Recent  examples are the combinatorial  computation of its Euler characteristic \cite{bodo}; 
the combinatorial aspects in  the study of the Nef cone of $\overline M_g$  (see in particular question 0.13 of \cite{GKM}). 
The importance of combinatorics is evident also in the study of spin curves \cite{CC} \cite{CCC}.
%

The main result relating the geometry of $\overline M_g$ to the combinatoric of stable curves is 
 is the stratification of $\overline M_g$ by topological type, 
which is governed by the {\em weighted dual graph} associated  to any stable curve 
 (Definition \ref{dwg}).
To any (multi)graph, 
it is naturally associated a group,
which we will call {\em complexity group} of the graph (Definition \ref{cg}).
In particular, given a stable (or, more generally, nodal) curve $C$ we call $\Delta_C$ the complexity group of the
dual graph of the curve.
This group has been extensively studied as an invariant of graphs, with applications to
 Physics, Chemistry, and many other areas,  and  it goes under many different names, 
 such as critical group \cite{biggschip} \cite{C-R}, determinant group \cite{tatiana}, Picard group \cite{tatiana}, Jacobian group \cite{B-N},  
 abelian sandpiles group \cite{C-E}.
From the point of view of geometry, the complexity  group was introduced in \cite{cap}, with the name of {\em degree class group},
in order to describe and handle the fibres of the compactification of the universal Picard variety $\overline P_{d,g}$ 
over $\overline M_g$ (also constructed in the same article). 
%
Moreover, this group arises naturally in the study of the 
compactified Jacobians of families of curves.
More precisely, let $R$ be a discrete valuation ring with fraction field $K$ and residue field ${\kapp}$ and denote by $B$ the spectrum of $R$. 
Consider a flat and proper regular curve $X\longrightarrow B$, such that the closed fiber $X_{\kappy}$
is geometrically irreducible.
Under some technical assumptions,  it is defined  the  group $\Phi$ of connected
components of the  the N\'eron model of the  Jacobian $J_K$ of the generic fiber $X_K$
(see \cite{BLR} sec. 9.6).
Notice that in this definition the closed fiber doesn't need to be nodal;
when $X_{\kappy}$ is nodal, $\Phi$ is precisely $\Delta_{X_{\kappy}}$.
Caporaso in  \cite{capneron} gave a geometric counterpart of this construction, showing that, under the assumption that $(d-g+1, 2g-2)=1$, there exists a space over $\overline M_g$ such that, for every regular family of stable curves over $B$, the N\'eron model of the Picard variety of degree $d$ of $X_K$ is obtained by base change via the moduli map $B\longrightarrow \overline M_g$. A yet another geometric interpretation of $\Phi$ is given by Chiodo in \cite{chiodo}, where the $r$-torsion points of $\Phi$ are described as N\'eron models of $r$-torsion line bundles on $ X_K$.

The relationship between the structure of this group and the structure of the graph has been studied, among others, by Lorenzini
in several papers (see for instance  \cite{lorgraphs} and \cite{lorjac}).
Other important references are \cite{biggschip},  \cite{tatiana} and \cite{C-R}.
We addressed this problem in \cite{BMS}, where in particular we constructed
a family of graphs with cyclic complexity group.
It seems that the question of finding a relation between the structure of the complexity group and the geometry of the curve 
is extremely difficult and intricate.

Rather that in the structure, we are interested here in the cardinality of the complexity group.
Kirkoff's Matrix Tree Theorem says that this 
integer is the number of spanning trees of the dual graph of the curve, usually called the {\em complexity}; 
this is the reason for our notation, also suggested by L. Caporaso in \cite{prag}. 
Moreover,  it can be calculated 
as the determinant of a certain  matrix (Theorem \ref{MTT}, Remarks \ref{spanningtrees} and \ref{bourbaki}).

We shall call complexity of a curve the corresponding complexity of its dual graph.
The complexity being an upper semicontinuous function over $\overline M_g$ (Lemma \ref{usc}), 
it defines a weak stratification on $\overline M_g$.
In Section \ref{pre} we investigate the relationship between this stratification and the (strong) one given by topological type.
This relation is also enlighted using the list of possible graphs for curves of genus $3$ given in Section \ref{lista}.

In particular, we are interested in the classes of curves with {\em maximal complexity}. Define the function
$\psi(g):=\mbox{max}\{ |\Delta_C|,\,\, C \mbox{ stable curve of genus } g\}$;
our guiding problems are the following.\\

1. Give a characterization of the curves $C$ such that $c(C)=\psi (g)$, or at least  with ``big'' complexity.\\

2. Find  bounds for $\psi$, depending only on $g$.\\

In Section \ref{MAX}, we give partial answers to problem 1, 
finding necessary conditions for curves to have maximal complexity.
We show that these curves are all graph curves in the sense of \cite{B-E}  i.e., they have simple trivalent dual graph 
$\Gamma$ with $b_1(\Gamma)=g$; moreover, they have no disconnecting nodes (Theorem \ref{maxconf} and its combinatorial version Theorem \ref{combinmaxconf}).
Using this reduction, we can apply some results of  the immense literature regarding  the complexity of 
regular graphs.

In Section \ref{ub} we give an answer to our second problem, 
 using  results of Biggs \cite{biggs}, McKay \cite{McK} and Chung-Yau \cite{C-Y}.
In particular, we obtain an asymptotically sharp upper bound  for $g\gg 0$.
By what observed in the beginning, these bounds limit the number of connected components of the 
N\'eron model of the generalized Jacobian  of stable curves. 
Note that this result is different from the one given by Lorenzini in \cite{lorner}: see Remark \ref{confrontoL}.
%
In Section \ref{lower}, we provide an example of families of graphs corresponding to stable curves with increasing genus,
and compute explicitly the complexity depending on the genus, thus giving an explicit lower bound.

In the last section we discuss the conjectural behavior of the curves with maximal complexity.
In particular, it seems that the graph of  such a curve should have maximal connectivity (which is 3 in this case).
Moreover, it seems that the girth has to be big (Conjecture \ref{connectivity} and \ref{girth}).
We prove in this section some partial results that seem to support the conjectures (Proposition \ref{bah} and Corollary \ref{g<10}).
Moreover, we prove a uniform necessary condition holding for any sequence of graphs of curves  with maximal complexity
(Theorem \ref{cicli}), using again a result of McKay.

Even though our point of view is irreparably geometrical, our main results are proven with (very simple) combinatorics methods, 
and  can be rewritten in purely combinatorial terms (see in particular \ref{combi}).
We would like to express the hope that our geometric approach does not discourage non-geometers from the reading of the article;
and in particular from considering the list of open question given in section \ref{congetture}.

\medskip

\noindent {\bf Acknowledgments}

We would like to express our gratitude to L. Caporaso for introducing us the problems we are dealing with in the present paper and for following this work very closely
and to C. Casagrande for the important corrections and suggestions she gave us.

We also thank A. Machi for transmitting us much encouragement after reading a preliminary version.

\section{Combinatorial invariants of nodal curves}\label{pre}

We work over an algebraically closed field $\kapp =\overline \kapp$.
A {\em nodal} curve is a reduced curve which has only ordinary double points as singularities.

\subsubsection*{The dual graph of a curve}
\begin{defi}\label{dwg}
To a nodal curve $C$ we can associate a graph $\Gamma_C$,
called the \textit{weighted dual graph}, 
given by a triple $(V,E,g)$, where 
$V$ is the set of vertices, $E$ the set of edges, and $g$ a function on the set  $V$ with non-negative integer values.
This triple is defined in the following way
\begin{itemize}
\item to each irreducible component $A$ corresponds a vertex\footnote{
Sometimes in the literature the vertices of a graph are denoted ``nodes''; of course we will never adopt this notation
which is extremely confusing in our context. A node will always be for us an ordinary double point of a curve!} 
$v_A$; 
\item to each node intersecting the components $A$ and $B$ (where $A$ and $B$ can coincide) 
corresponds an edge connecting the vertices $v_A$ and $v_B$;
\item $g: V\longrightarrow \Z_{\geq 0}$ is 
 the function that associates to any vertex $v$ the geometric genus of the corresponding component of $C$.
\end{itemize}
\end{defi}
Call $\gamma$ the number of irreducible components of $C$ and 
$\delta$ the number of nodes.
Thus $\Gamma_C$ has $\gamma$ vertices, $\delta$ edges, 
and among the edges there is a loop 
for every node lying on a single irreducible component of $C$.
The weighted graph encodes all the topological 
information about the curve.
Of course, conversely, any weighted graph can be realized as the dual graph of a nodal curve.




It is important to stress that dual graphs of nodal curves can have more than one edge connecting two nodes, and can also have loops.
These are usually called {\em multigraphs}. 
In this paper, by graph we will always mean multigraph, while a graph without loops and multiple edges will be called {\em simple}.

Given a graph $\Gamma$ with $\delta$ edges, $\gamma$ vertices and $c$ connected components, its first Betti number is
$b_1(\Gamma):= \delta-\gamma+c$; it corresponds to the number of independent cycles on $\Gamma$.

Let $\{C_1,\ldots,C_{\gamma}\}$ be the set of  irreducible components of a nodal curve $C$ 
with $\delta$ nodes and $c$ connected components.
Recall that the arithmetic genus of $C$ can be computed by the following formula
\begin{equation}\label{genere}
g(C)=\sum_{i=1}^{\gamma}g(C_i)+\delta-\gamma+c=\sum_{v\in V(\Gamma_C)}g(v)+b_1(\Gamma_C).
\end{equation}
By analogy, for any weighted graph $\Gamma$ given by the triple  $(V,E,g)$, we will call 
$g(\Gamma):=\sum_{v\in V}g(v)+b_1(\Gamma)$ the {\em (arithmetic) genus of the graph}.


\subsubsection*{Complexity group} 

Let us consider a connected nodal curve $C$.
Let $\left\{C_i\right\}_{i=1,...,\gamma}$ be the irreducible components of $C$. 
Define 
$$
\begin{matrix}k_{ij}:= & \left\{\begin{array}{l}\;\;\, \sharp(C_i\cap C_j) \, \mbox{ if }i\not= j\\  \\
-\sharp(C_i\cap \overline{C\setminus C_i}) \, \mbox{ if } i=j\\\end{array}\right.\end{matrix}
$$
As $C_i\cap \overline{C\setminus C_i}=\bigcup_{j\not= i} C_i\cap C_j$, we have that, for fixed $i$, $\sum_j k_{ij}=0$.
For every $i$ set 
$$\underline{c}_i:=(k_{1i},\ldots,k_{\gamma i}) \in \mathbb Z^\gamma.$$ 
Call $Z:=\{\underline{z}\in\mathbb{Z}^\gamma : |\underline z|=0\}$.
As observed before, $\underline c_i\in Z$. Let us call $\Lambda_C$ the sublattice of $Z$ spanned by 
$\{\underline c_1,\ldots,\underline c_\gamma\}$. In fact, $\Lambda_C$ is a lattice in  $Z$ (it
has rank $\gamma -1$) as we will show in a moment (see \cite{capneron} for a geometric proof of this fact).

\begin{defi}\label{cg}
The complexity group of $C$ is the finite abelian group $\Delta_C:=Z/\Lambda_C$.
\end{defi}
It is important to notice that this group depends {\em only} on the dual (non-weighted) graph
of the curve; clearly it is defined  for any connected\footnote{The definition could be easily extended to 
non connected graphs/curves.}  graph.
As noted in the Introduction, this group arises in many  contexts where graphs are used,
and it is known with many other names.

Let $M$ be the $\gamma\times\gamma$ matrix whose columns are the $\underline c_i$'s. 
We will call $M$ the {\em intersection matrix}, the name clearly deriving from its geometrical meaning. 
However, in literature, $M$ is known as the {\em (combinatorial) Laplacian} matrix (cf. e.g. \cite{bollobas} and \cite{lorfinite}). 
It is obtained from the so-called adjacency matrix of $\Gamma_C$ subtracting the vertex degrees on the diagonal.


\subsubsection*{Complexity of a graph}

A tree is a connected graph $G$ with $b_1(G)=0$.
Let $\Gamma$ be a graph. A {\it spanning tree} of $\Gamma$ is a subgraph of $\Gamma$ which is a tree having 
the same vertices as $\Gamma$. 
\begin{defi}
The complexity of $\Gamma$, indicated by the symbol $c(\Gamma)$, is the number of spanning trees contained 
in $\Gamma$ (see e.g. \cite{biggs}, sec 6, \cite{berge}, cap.3 $\natural$ 5,  \cite{west}, sec 2.2).
\end{defi}

Observe that $c(\Gamma)=0$ if and only if $\Gamma$ is not connected, 
and that if $\Gamma$ is a connected tree  $c(\Gamma)=1$.
For the complexity of the dual graph associated to a curve $C$, we will often use the symbol $c(C)$, instead of
$c(\Gamma_C)$.
The following theorem, known as Kirkoff's Matrix Tree Theorem, will be a key ingredient for our work. 
There are at least three different proofs of this result; see \cite{BMS} for a proof and for the references.

\begin{teo}{\upshape (Matrix Tree Theorem)}\label{MTT}
Let $s, t \in \{1,\ldots \gamma\}$.
Using the above notations, if $M^\star$ is obtained from $M$ by deleting the $t$-th column and the $s$-th row, then 
$$c(\Gamma)=(-1)^{s+t+\gamma -1}\mbox{det}(M^\star).$$
\end{teo}

In particular, the Matrix Tree Theorem assures that, in the case of the 
dual graph of a {\em connected} curve $C$, the matrix  $M$ has rank $\gamma-1$
i.e., $\Lambda_C$ is indeed a lattice.

\begin{rem}\label{spanningtrees} \upshape{
For $r\in\{1,\ldots,\gamma\}$, consider the isomorphism $\alpha_r:Z\stackrel{\sim}\longrightarrow\mathbb{Z}^{\gamma - 1}$
which consists of deleting the $r$-th component.
The group $\Delta_C$ is the quotient of $\mathbb{Z}^{\gamma - 1}$ by the lattice generated by 
$$\underline c_i^\prime:=(k_{1i}
,\ldots,\widehat{k_{ri}},\ldots,k_{\gamma i}).$$
Observe that again $\sum_i \underline c_i^\prime = \underline 0 \in \mathbb{Z}^{\gamma -1}$.
Therefore $\Delta_C$ is presented by the matrix $M^\star$ obtained from $M$ by deleting a column and the $r$-th row.
Hence, we can compute its cardinality via Theorem \ref{MTT}
$$c(\Gamma_C)=|\Delta_C|=|\mbox{det}(M^\star)|.$$
So, we can conclude that {\em the cardinality of the complexity group of a curve $C$ is the complexity of the dual graph $\Gamma_C$.}}
\end{rem}

\begin{rem}\label{bourbaki}\upshape{
By the diagonalization theorem of integer matrices (i.e. the structure theorem for finite abelian groups, 
see \cite{ar}), $M$ is equivalent over $\mathbb Z$ to a diagonal matrix
$diag (d_1, d_2,\ldots, d_\gamma)$, where the $d_i's$ are the invariant factors, i.e. the 
non-negative integers  obtained in the following way:
Let $D_i$ be the greatest common divisor of the $i\times i$ minors of $M$. 
Then  $d_i=D_i/D_{i-1}$
(cf. also \cite{lorarithmetical}, Theorem 1.5). 
Note that $d_\gamma=0$ and $d_i>0$ for $i\not = \gamma$, so 
$\Delta_C=\oplus _{i=0}^{\gamma-1}\mathbb Z/d_i \mathbb Z $.
In particular, $ |\Delta_C|$ is equal to the greatest common divisor of all the $(\gamma-1)\times(\gamma-1)$ 
minors of $M$. 

On the other hand, if we diagonalize $M$ over the real numbers, we get real eigenvalues
$0=\lambda_1>\lambda_2\geq \ldots\geq \lambda_\gamma$. 
These eigenvalues are deeply studied in Combinatorics.
Note that in general the $\lambda_i$'s have no relation with  the invariant factors, not even if they happen to be integers; 
a nice counterexample can be found in Section 9.2 of \cite{tatiana}.
An attempt to construct the invariant factors from the $\lambda_i's$, for a particular class of graphs, can be find in \cite{C-R}.

Let us note that in particular $ | \Delta_C|=\gamma^{-1}\lambda_2\lambda_3\ldots\lambda_\gamma$ 
(cf. \cite{biggs} cor.6.5).
}
\end{rem}


\subsubsection*{Stable curves and their graphs}

\begin{defi}\label{stable}
A curve $C$ of genus $g\geq 2$ over $\kapp$ is stable (resp. semistable) if it is  nodal, connected and  such that if 
$D\subset C$ is a smooth rational component, then $|D\cap\overline{C\setminus D}|\geq 3$ (resp. $\geq 2$).
\end{defi}

The moduli space of stable curves of genus $g$, $\overline M_g$, is a projective variety of dimension $3g-3$, the so-called Deligne-Mumford compactification of the moduli space of smooth curves of genus $g$, $M_g$.
%
The theory of stable curves was first introduced by A. Mayer and D. Mumford and was first developed in \cite{DM} in order to prove the irreducibility of $M_g$ in any characteristic. In that paper, the authors prove the main properties of stable curves used later by D. Gieseker in \cite{gie} to establish the existence of $\overline M_g$.

\medskip

Given a graph $\Gamma$ and a vertex $v$, the degree $d(v)$ of $v$  is the number of half edges touching $v$.
The combinatorial version of the definition of a stable curve is the following.

\begin{defi}\label{cstable}
A weighted graph  $\Gamma$ of arithmetic genus $g\geq 2$ is stable if
\begin{equation}\label{stability}
2g(v)-2+d(v)>0 \,\,\,\,\,\,\mbox{ for any }\,\,\,\,\,v\in V.
\end{equation}
$\Gamma$ is said to be semistable if the inequality above holds with $\geq$.
\end{defi}


\subsubsection*{Topological stratification of $\overline M_g$}

There is a natural stratification on $\overline M_g$ given by topological type.
Each stratum of codimension $k$ is the subset of $\overline M_g$ 
consisting of classes of  stable curves  with weighted graphs of genus $g$  having exactly $k$ edges.
In particular, the stratum of codimension $0$ is the open set $M_g\subset \overline M_g$ of smooth curves.
The strata of codimension $1$ are $[g/2]+1$, and the corresponding graphs are of the following types.

\bigskip

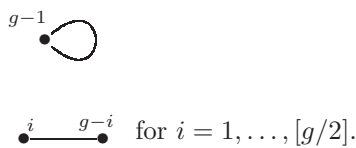
\centerline{
\xymatrix@=.2pc{
*{\bullet} \ar@{-}^<{g-1} @(ur,dr)\\{}
}
\,\,\,\,\,\,\,\,\,\,\,\,\,\,\,\,\,\,\,\,\,\,\,\,\,\,\,\,\,\,\,\,\,\,\,\,\,\,\,\,\,\,\,\,\,\,\,\,\,\,\,\,\,\,\,\,\,\,\,\,\,\,\,\,\,\,\,\,\,\,\,\,
}

\bigskip
\centerline{
\xymatrix@R=.1pc{
*{\bullet} \ar@{-}^<{i}^>{g-i}[r] & *{\bullet}
\\{}
}
\,  for $i=1,\ldots, [g/2]$.
}

\bigskip

\noindent The closure of these strata are effective divisors on $\overline M_g$, the so-called boundary divisors,
usually denoted by $\Delta_i$, for  $i=0,\ldots, [g/2]$, according to the preceding list.
For instance, $\Delta_0$ has as generic point  corresponding to an irreducible curve of geometric genus $g-1$ 
with exactly one node.
It can be described as the locus 
of curves having at least one non-disconnecting node.
The boundary of $\overline M_g$ is $\partial \overline M_g=\overline M_g\setminus M_g=\cup_{i=0}^{[g/2]}\Delta_i$.

The union of the strata of  codimension $k$ is the variety  of isomorphism classes of stable curves having exactly $k$ nodes; 
its closure is the locus of curves with at least $k$ nodes.

\subsubsection*{Degeneration of nodal curves}

The following well-known result describes the possible transformations on a dual graph corresponding to the degenerations of a nodal curve.

\begin{prop}\label{degenerations}
Let $\Gamma=(V,E,g)$ be a weighted graph. 
Let  $C$ be a nodal  curve having $\Gamma$ as dual weighted  graph.
The dual graphs of the possible nodal curves obtained by degenerations of $C$ 
are obtained from $\Gamma$ via a sequence of the following
transformations:
\begin{itemize}
\item[(I)] given a vertex $v$ such that $g(v)\geq 1$ add a loop on $v$ and decrease its weight by 1;
\item[(II)] given a vertex $v$, and given two nonnegative integers $a$ and $b$ such that $a+b=g(v)$, substitute it with 
two vertices $v_a$ and $v_b$ with weights respectively $a$ and $b$  and one edge $l$ connecting them.
\end{itemize}
\end{prop}
\noindent The figure below illustrates operation (II).

\vspace{0.5cm}
\centerline{
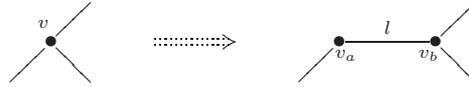
\begin{tabular}{cc}
{}
&
\xymatrix@=1pc{
& &\ar@{-}[dl]_>v &&&& &&&&\ar@{-}[dl]\\
& *{\bullet} && \ar@{:>}[rr] && & &*{\bullet} \ar@{-}[rr]^l_<{v_a}_>{v_b}&& *{\bullet}& \\
\ar@{-}[ur]& &\ar@{-}[ul] &&&& \ar@{-}[ur] &&&&\ar@{-}[ul]
}
\end{tabular}}
\vspace{0.5cm}

\noindent Of course, there are many ways to perform operation (II) (see also Definition \ref{operazione} below).
Notice that this is the opposite operation of contracting an edge, in the sense that if we contract the 
new edge $l$ we get the original graph.

Note also that if our given curve $C$ is stable, degenerations of $C$ obtained by performing operation (I) are still stable, while if we perform operation (II) this is the case only if
$2a-2+d(v_a)\geq 1$ and $2b-2+d(v_b)\geq 1$.

\subsubsection*{Polygonal curves}

We present here an elementary combinatorial proof of the following well-known fact for stable curves. 
A geometric proof can be obtained using the above results about the topological stratification of $\overline M_g$.
\begin{lem} \label{lemcomb}
Let $C$ be a stable curve of genus $g\ge 2$. 
Then
\begin{enumerate}
\item $C$ has at most $3g-3$ nodes and $2g-2$ irreducible components.
\item Assume that $C$ has $3g-3$ nodes. Then $C$ has $2g-2$ components $C_1,\ldots,C_{2g -2}$ and, 
if $\nu_i:C_i^\nu\longrightarrow C_i$ is the normalization of $C_i$, then 
$C_i^\nu\simeq \mathbb{P}^1$ and $|\nu_i^{-1}(C_i\cap C_{sing})|=3$ for all $i$.
\end{enumerate}
\end{lem}
\begin{proof}
Let $\Gamma=(V,E,g)$ be the weighted graph associated to $C$.
Note that $\sum_{v\in V}d(v)=2\delta$, hence we can rewrite formula (\ref{genere}) as
\begin{equation}\label{genere2}
g=\sum_{v\in V}\left(g(v)+\frac{d(v)}{2}\right)-\gamma +1.
\end{equation}
By the connectedness of $C$, $d(v)\geq 1$ for any $v$. 
Using also the stability condition, we have that 
$
g\geq \sum_{v\in V}\frac{3}{2}-\gamma +1=\frac{\gamma}{2}+1,
$
hence, $\gamma \leq 2g-2$. Now, using (\ref{genere}) again, we have that $\delta \leq 3g-3$, and (1) follows.

Now, if $\delta =3g-3$, again by (\ref{genere}), we see that necessarily $\gamma =2g-2$ and 
$g(v)=0$ for any $v\in V$.
Moreover, using again (\ref{genere2}), we see that $d(v)=3$ for any vertex $v$, so also (2) is proved.
\end{proof}

Stable curves with $3g-3$ nodes (and $2g-2$ components), 
are the $0$-strata of the topological stratification of $\overline M_g$; 
they are rigid, in the sense that 
any deformation of such curves in a family of stable curves must have necessarily at least one node smoothed. 
Indeed, from  Proposition \ref{degenerations}  we see that both operations (I) and (II) cannot be performed, 
i.e. there is no possible further degeneration for such curves.
We will call such curves {\em polygonal curves}
(also called large limit curves, see \cite{tyurin}).
Polygonal curves with simple graph are called graph curves (\cite{B-E}).

We refer to \cite{B-E} for an ample discussion on the importance and the properties of these curves.
Let us just make a couple of remarks. 
The associated dual graphs are trivalent (multi)graphs such that the weight function $g$ is $0$ on every vertex.
Due to the fact that $\mathbb P^1$ with $3$ fixed points has no non-trivial automorphisms, 
the automorphism group of such curves coincides with the automorphism group of the graph.
Artamkin in \cite{artamkingen} has given a recursive rule that computes the number
$$
\sum_{\begin{array}{c}\Gamma \,\,\,trivalent,\\
b_1(\Gamma)=g
\end{array}}
\frac{1}{|Aut(\Gamma)|} \mbox{ ,}
$$
which is the ``stacky top self-intersection'' of the boundary divisor $\partial \overline M_g$.

\subsubsection*{Stratification by complexity of $\overline M_g$}

We can define a  function $c\colon \overline M_g\longrightarrow \mathbb Z_{\geq 0}$ 
associating to every $[C]\in \overline M_g$ its complexity $c(C)$.

Of course, $c$ is bounded from above, i.e., there is a bound for the complexity of a stable curve of given genus;
indeed, by Lemma \ref{lemcomb}, there is only a finite number of possible graphs.
In section \ref{ub} we provide upper bounds, which are asymptotically sharp for $g\gg 0$.
\begin{rem}\upshape{
Clearly, this wouldn't make sense for nodal curves (not even for semistable ones). 
Indeed, blowing up a node an arbitrary number of times doesn't change the genus of the curve, 
but it can increase arbitrarily the complexity of a nodal curve (see for instance \cite{BMS} Proposition 3.3.).}
\end{rem}

\begin{lem}\label{usc}
The function $c\colon \overline M_g\longrightarrow \mathbb Z_{\geq 0}$  is upper semi-continuous.
\end{lem}
\begin{proof}
Follows from Proposition \ref{increasing} below.
\end{proof}
This result implies in particular that we could define a stratification ``by complexity'' of $\overline M_g$.
This is clearly a much rougher stratification than the one by topological type; 
{\em the set of curves in $\overline{M}_g$ with given complexity is a (maybe empty) union of components of different codimension strata 
of the topological type stratification, of different codimension} (see the case of genus $3$ in the next section).


For instance, the set $\overline M_g^{c=1}$ of curves with complexity one is the set of curves whose dual graph is a tree with loops.
It contains the curves of compact type $\overline M_g^{ct}=\overline M_g\setminus \Delta_0$, but also the interior of $\Delta_0$,
and other strata of bigger codimension; it also contains $0$-strata, i.e. isolated points (see next section).
The complement $\overline M_g\setminus \overline M_g^{c=1}$ of curves with complexity greater than one is a closed 
subset of codimension $2$.


\section{List of graphs for $\overline M_3$}\label{lista}

In this section, we list all the possible weighted graphs for stable curves of genus 3, 
as well as their complexity and their complexity group. 
This list (with one error, now corrected, which was kindly pointed out to us by A. Chiodo) appeared in \cite{BMS}.
We will use $\mathbb{Z}_n$ to denote the quotient group $\mathbb{Z}/n\mathbb{Z}$. 
The graphs are ordered by increasing the number of nodes (i.e. according to the codimension of corresponding 
stratum of the stratification by topological type).
In the graphs we will indicate the weight of each vertex only if it is not zero.
Recall that if $C$ is a stable curve of genus 3, then $C$ has at most $6$ nodes and $4$ components.

\begin{equation*}
\begin{tabular}{||c|c|c|c|c||}
\hline \hline
Graph configuration & Nodes & Components & Complexity & DCG\\
\hline
\xymatrix@=.2pc{
{}\\
*{\bullet}\\{}
\scriptstyle{3}}
&
\xymatrix@=.001pc{
\\
0}
&
\xymatrix@=.001pc{
\\
1}
&
\xymatrix@=.001pc{
\\
1}
&
\xymatrix@=.001pc{
\\
0}\\

\xymatrix{
*{2 \,\bullet} \ar@{-}@(ur,dr) & {}
}
&
1
&
1
&
1
&
0 \\
&&&&\\

\xymatrix{
*{2\, \bullet} \ar@{-}[r] & *{\bullet \, 1}
}
&
1&2&
1
&
0\\
&&&&\\
\xymatrix{
*{\bullet} \ar@{-}_<{1}@(ul,dl) \ar@{-}@(dr,ur)
}
&
2&1&
1
&
0
\\
&&&&\\

\xymatrix{
*{\bullet} \ar@{-}@(ul,dl) \ar@{-}[r]^<{1}^>{1}& *{\bullet}
}
&
2&2&
1
&
0\\
&&&&\\
\xymatrix{
*{\bullet} \ar@{-}@(ul,dl) \ar@{-}[r]^>{2}& *{\bullet}
}
&
2&2&
1
&
0
\\
&&&&\\
\xymatrix{
*{\bullet} \ar @{-} @/_/[r]_>{1}  & *{\bullet} \ar@{-} @/_/[l]_>{1}
}
&
2&2&
2
&
$\mathbb{Z}_2$\\
&&&&\\
\xymatrix@R=1pc{
*{\bullet} \ar@{-}[r]^<{1} & *{\bullet}\ar@{-}[r]^<{1}^>{1} & *{\bullet}\\
{}
}
&
2&3&
1
&
0\\
\hline
\hline
\end{tabular}
\end{equation*}

\begin{equation*}
\begin{tabular}{||c|c|c|c|c||}
\hline \hline
Graph configuration & Nodes & Components & Complexity & DCG\\
\hline 
{}&&&&\\
{}&&&&\\
\xymatrix@R=1pc{
*{\bullet} \ar@{-}@(ul,dl) \ar@{-}@(ur,dr) \ar@{-}@(ur,ul)\\
{}
}
&
3&1
&
1&
0\\
\xymatrix@R=1pc{
*{\bullet} \ar@{-}@(ur,ul) \ar@{-}@(dr,dl) \ar@{-}[r]^>{1}& *{\bullet}\\
{}
}
&
3&2&
1
&
0\\
\xymatrix{
*{\bullet} \ar@{-}@(ul,dl) \ar@{-}[r]_<{1} & *{\bullet} \ar@{-} @(ur,dr)
}
&
3&2&
1
&
0\\
&&&&\\
\xymatrix{
*{\bullet} \ar @{-}@(ul,dl) \ar @{-}@/_/[r] & *{\bullet} \ar @{-}_<{1}@/_/[l]
}
&
3&2&
2
&
$\mathbb Z_2$\\
&&&&\\
\xymatrix{
*{\bullet} \ar @{-} @/_/[r] \ar@{-}[r] & *{\bullet} \ar@{-}_>{1} @/_/[l]
}
&
3&2&
3
&
$\mathbb Z_3$\\
&&&&\\
\xymatrix@R=.5pc{
*{\bullet} \ar@{-}@(ul,dl) \ar@{-}[r] & *{\bullet}\ar@{-}[r]^<{1}^>{1} & *{\bullet}\\
{}
}
&
3&3
&
1
&
0\\
\xymatrix@R=.2pc{
*{\bullet} \ar@{-}[r]^<{1} & *{\bullet}\ar@{-}@(ul,ur) \ar@{-}[r]^>{1} & *{\bullet}\\
{}
}
&
3&3&
1&
0\\
\xymatrix{
 *{\bullet} \ar@{-}^<{1}[r]& *{\bullet} \ar @{-} @/_/[r] & *{\bullet} \ar@{-}_<{1}@/_/[l]
}
&
3&3&
2
&
$\mathbb Z_2$\\
&&&&\\
\xymatrix@R=1pc{
*{\bullet} \ar@{-}@(ur,ul) \ar@{-}@(dr,dl) \ar@{-}[r]& *{\bullet} \ar@{-}@(ur,dr)\\
{}
}
&
4&2&
1&
0\\
\xymatrix{
*{\bullet} \ar @{-}@(ul,dl) \ar @{-}@/_/[r] & *{\bullet} \ar @{-}@/_/[l] \ar @{-}@(ur,dr)
}
&
4&2&
2
&
$\mathbb Z_2$\\
&&&&\\
\xymatrix{
*{\bullet} \ar@{-}@(ul,dl) \ar @{-} @/_/[r] \ar@{-}[r] & *{\bullet} \ar@{-}@/_/[l]
}
&
4&2&
3
&
$\mathbb Z_3$\\
&&&&\\
\xymatrix{
*{\bullet}  \ar @{-} @/_.3pc/[r] \ar @{-}@/_.8pc/[r]/& *{\bullet} \ar@{-}@/_.3pc/[l] \ar@{-}@/_.8pc/[l]
}
&
4&2&
4
&
$\mathbb Z_4$\\
&&&&\\
\xymatrix@R=.5pc{
*{\bullet} \ar@{-}@(ul,dl) \ar@{-}[r] & *{\bullet}\ar@{-}[r]^<{1} & *{\bullet}\ar@{-}@(ur,dr)\\
{}
}
&
4&3
&
1
&
0\\
\xymatrix@R=.2pc{
*{\bullet} \ar@{-}[r] \ar@{-}@(ul,dl) & *{\bullet}\ar@{-}@(ul,ur) \ar@{-}[r]^>{1} & *{\bullet}\\
{}
}
&
4&3&
1&
0\\
\xymatrix@R=.1pc{
{} & *{\bullet} \ar@{-}@(ul,dl) \ar@{-}[r]& *{\bullet} \ar @{-} @/_/[r] & *{\bullet} \ar@{-}_<{1}@/_/[l]\\
{}
}
&4&3
&
2
&
$\mathbb Z_2$\\
\xymatrix@=1.2pc{
& *{\bullet} \ar@{-}[dl] \ar@{-}^<{1}[dr]& \\
*{\bullet} \ar@{-}[rr] \ar@{-}@/_/[rr]& & *{\bullet}
}
&
\xymatrix@=.3pc{
\\
4}
&
\xymatrix@=.3pc{
\\
3}
&
\xymatrix@=.3pc{
\\
5\\
{}
}
&
\xymatrix@=.3pc{
\\
{\mathbb Z_5}
}
\\
\xymatrix@R=.8pc{
{} & *{\bullet} \ar@{-}^<{1}[r]& *{\bullet} \ar@{-}[r] \ar @{-} @/_/[r] & *{\bullet} \ar@{-}@/_/[l]&
}
&
4&3&
3
&
$\mathbb Z_3$
\\
\xymatrix@=.5pc{
*{\bullet}\ar@{-}^<{1}[dr]&&&\\
&*{\bullet}\ar@{-}[rr]&&*{\bullet}\ar@{-}@(ur,dr)\\
*{\bullet} \ar@{-}_<{1}[ur]&&
}&
\xymatrix@=.3pc{
\\
4
\\
{}
\\
{}}
&\xymatrix@=.3pc{
\\
4}
&
\xymatrix@=.3pc{
\\
1}
&
\xymatrix@=.3pc{
\\
0}
\\
\xymatrix@R=.3pc{
*{\bullet} \ar@{-}[r] \ar@{-}@(ul,dl) & *{\bullet}\ar@{-}@(ul,ur) \ar@{-}[r] & *{\bullet}
\ar@{-}@(ur,dr)\\
{}
}
&
5&3&
1&
0\\
\xymatrix@R=.3pc{
{}&&\\
*{\bullet} \ar@{-}[r] \ar@{-}@(ul,dl) & *{\bullet} \ar@{-}@/^/[r] \ar@{-}@/_/[r] & *{\bullet}
\ar@{-}@(ur,dr)\\
{}
}
&
\xymatrix@R=.3pc{
{}\\
5}
&
\xymatrix@R=.3pc{
{}\\
3}
&
\xymatrix@R=.3pc{
{}\\
2}
&
\xymatrix@R=.3pc{
{}\\
\mathbb{Z}_2}\\

\hline \hline
\end{tabular}
\end{equation*}

\begin{equation*}
\begin{tabular}{||c|c|c|c|c||}
\hline \hline
Graph configuration & Nodes & Components & Complexity & DCG\\
\hline 
{}&&&&\\
\xymatrix@R=.2pc{
{}\\
{} & *{\bullet} \ar@{-}[r] \ar@{-}@(ul,dl)& *{\bullet} \ar@{-}[r] \ar @{-} @/_/[r] & *{\bullet} \ar@{-}@/_/[l]&\\&&&&
}
&
\xymatrix@R=.2pc{
{}\\
5}
&
\xymatrix@R=.2pc{
{}\\
3}
&
\xymatrix@R=.2pc{
{}\\
3}
&
\xymatrix@R=.2pc{
{}\\
\mathbb Z_3}
\\
\xymatrix@=1.2pc{
& *{\bullet} \ar@{-}[dl] \ar@{-}@/_/[dl] \ar@{-}@/^/[dr] \ar@{-}[dr]& \\
*{\bullet} \ar@{-}[rr] & & *{\bullet}\\{}
}
&
\xymatrix@=.3pc{
\\
5}
&
\xymatrix@=.3pc{
\\
3}
&
\xymatrix@=.3pc{
\\8
\\
{}
}
&
\xymatrix@=.3pc{
\\
{\mathbb Z_8}}\\
\xymatrix@=1.2pc{
& *{\bullet} \ar@{-}[dl] \ar@{-}[dr] \ar@{-}@(ul,ur)& \\
*{\bullet} \ar@{-}[rr] \ar@{-}@/_/[rr]& & *{\bullet}
}
&
\xymatrix@=.5pc{
\\5}
&
\xymatrix@=.5pc{
\\3
\\
{}
}
&
\xymatrix@=.5pc{
\\5
\\
{}
}
&
\xymatrix@=.5pc{
\\
{\mathbb Z_5}}
\\
\xymatrix@=1.2pc{
&*{\bullet} \ar@{-}^<{1}[d]&\\
& *{\bullet} \ar@{-}[dl] \ar@{-}[dr] & \\
*{\bullet} \ar@{-}[rr] \ar@{-}@/_/[rr]& & *{\bullet}
}
&\xymatrix@=.5pc{
{}\\
{}\\
5
}
&
\xymatrix@=.5pc{
{}\\
{}\\
4}
&
\xymatrix@=.5pc{
{}\\
{}\\
5
\\
{}\\{}
}
&
\xymatrix@=.5pc{
{}\\
{}\\
{\mathbb Z_5}}\\
\xymatrix@R=.2pc{
{} & *{\bullet} \ar@{-}[r] \ar@{-}@(ul,dl)& *{\bullet} \ar @{-} @/_/[r] & *{\bullet} \ar@{-}@/_/[l] \ar@{-}^>{1}[r]& *{\bullet}
\\&&&&
}
&5&4&
2
&
$\mathbb Z_2$
\\
\xymatrix@=.5pc{
*{\bullet}\ar@{-}^<{1}[ddd]&&&\\
&&*{\bullet}\ar@{-}@(ur,dr)&\\
&&&\\
*{\bullet}\ar@{-}[rrr]\ar@{-}[rruu]&& &*{\bullet}\ar@{-}@(ur,dr)\\
{}
}
&\xymatrix{
\\5}
&\xymatrix{
\\4}
&
\xymatrix{
\\1}
&
\xymatrix{
\\
0}
\\
\xymatrix@R=.2pc{
*{\bullet} \ar@{-}[r] \ar@{-}@(ul,dl)& *{\bullet} \ar @{-} @/_/[r] & *{\bullet} \ar@{-}@/_/[l] \ar@{-}[r]& *{\bullet} \ar@{-}@(ur,dr)
\\{}
}
&6&4&
2
&
$\mathbb Z_2$
\\
\xymatrix@=1pc{
*{\bullet} \ar@{-}[dd]\ar@{-}@/_.5pc/[dd]\ar@{-}[drr]&&&\\
&&*{\bullet}\ar@{-}[r]&*{\bullet} \ar@{-}@(ur,dr)\\
*{\bullet} \ar@{-}[urr]&&&
}
&
\xymatrix@=.6pc{
\\6}
&\xymatrix@=.6pc{
\\4}
&
\xymatrix@=.6pc{
{}\\
5
\\
{}\\{}
}
&
\xymatrix@=.6pc{
{}\\
{\mathbb Z_5}}\\

\xymatrix@=.5pc{
&*{\bullet}\ar@{-}[ddd] \ar@{-}@(ur,ul)&\\
*{\bullet} \ar@{-}@(ul,dl)&&*{\bullet}\ar@{-}@(ur,dr)&\\
&&&\\
&*{\bullet}\ar@{-}[ruu]\ar@{-}[luu] & \\
{}
}
&\xymatrix@=.5pc{
\\6}
&\xymatrix@=.5pc{
\\4}
&
\xymatrix@=.5pc{
\\1}
&
\xymatrix@=.5pc{
\\
0}
\\
\xymatrix{
*{\bullet} \ar@{-}[r] \ar@{-}[d] \ar@{-}@/_/[d] & *{\bullet} \ar@{-}[d] \ar@{-}@/^/[d]\\
*{\bullet} \ar@{-}[r] & *{\bullet}
}
&\xymatrix@=.5pc{
\\6}
&\xymatrix@=.5pc{
\\4}
&
\xymatrix@=.5pc{
\\12}
&
\xymatrix@=.5pc{
\\
\mathbb Z_{12}
}
\\
\xymatrix@=1.2pc{
&*{\bullet}\ar@{-}[dd] \ar@{-}[dl] \ar@{-}[dr] &\\
*{\bullet} \ar@{-}[rr] |!{[r]}\hole & &*{\bullet}\\
&*{\bullet}\ar@{-}[ru]\ar@{-}[lu] & 
}
&\xymatrix@=.5pc{
\\6\\
{}}
&\xymatrix@=.5pc{
\\4}
&
\xymatrix@=.5pc{
\\16}
&
\xymatrix@=.5pc{
\\
\mathbb Z_4\times\mathbb Z_4}
\\

\hline \hline
\end{tabular}
\end{equation*}

\bigskip

\bigskip

\bigskip

\noindent Note that the set of curves of genus $3$ with complexity $1$  contains at least one stratum 
of the topological stratification of any codimension.
The set of curves with complexity $2$ contains at least one stratum among the ones of codimension
greater or equal to $2$.
The set of curves with complexity $6$ is empty.


\section{Stable curves with maximal complexity}\label{MAX}

We will now focus our interest in curves with maximal complexity.
We can see that in case $g=3$ above, the curves reaching the maximal complexity are polygonal curves. 
However, as already observed, there are also polygonal curves with small, even $0$, complexity.
In what follows we prove that, indeed, any curve with maximal complexity is necessarily a polygonal curve
and, moreover, it is a graph curve without disconnecting nodes.
In Section \ref{congetture} we shall make some remarks and conjectures on sufficient conditions.

We shall need the following well-known  result, whose proof is elementary.

\begin{prop}\label{comb}
Let $\Gamma$ be a graph. If $e$ is an edge of $\Gamma$ which is not a loop,
call $\Gamma -e$ the graph obtained from $\Gamma$ by removing $e$, 
and $\Gamma\cdot e$ the one obtained by contracting $e$. 
Then, we have the following relation between the complexities of these three graphs:
\begin{equation}\label{fond}
c(\Gamma)=c(\Gamma - e)+c(\Gamma\cdot e).
\end{equation}
\end{prop}
\noindent See \cite{BMS} for a geometric interpretation and discussion of this proposition.\\


The main result of this section is the following

\begin{teo} \label{maxconf}
Let $C$ be a stable curve of genus $g\geq 3$ with maximal complexity. 
Then, $C$ is a curve without disconnecting nodes and with trivalent simple dual graph.
In particular, $C$ is a graph curve.
\end{teo}

From a geometrical point of view, this result implies in particular that the curves with maximal complexity 
lie on those $0$-strata of the topological stratification
which are contained in  $\partial \overline M_g\setminus \cup_{i=1}^{[g/2]} \Delta_i$.

See the end of this section for a purely combinatorial statement.

\bigskip

Our strategy is the following: we consider a generic stable curve $C$ of genus $g$ and its weighted dual graph, $\Gamma_C$;
then, if $\Gamma_C$ is not as stated in Theorem \ref{maxconf}, we modify the graph obtaining a new weighted graph $\Gamma'$, which is the dual graph of another stable curve 
of genus $g$ with the desired properties, and we prove that $c(\Gamma_C)<c(\Gamma')$.
So, we need to perform operations on the graph increasing strictly 
the number of spanning trees.
A. Kelmans made a deep study  of operations increasing the complexity of a graph, 
even though from a different point of view; see for instance \cite{Kelmans}.

Recall the two operations associated to the degeneration of a curve described in Proposition \ref{degenerations}.
We now study when do these operations increase the complexity of the graph. 

\begin{defi}\label{operazione}
We shall distinguish two different kinds of operation (II):
\begin{itemize}
\item[(II)a] if the new edge $l$ is a disconnecting edge;
\item[(II)b] otherwise, i.e. if there is  a cycle $\mathcal C$ that contains $v$ such that 
$\mathcal C\cup l$ is a cycle for the new graph.
\end{itemize}
\end{defi}

\begin{prop}\label{increasing}
Let $R$ be a discrete valuation ring with fraction field $K$ and residue field { $\kapp$}, 
and $f\colon X \longrightarrow \spc R$ a family of nodal curves.
Let $\Gamma_K$, resp. $\Gamma_{\small\kapp}$, be the weighted graph associated to the generic fiber $X_K$, 
resp. to the closed fiber $X_{\small\kapp}$. 

Then, the complexity function $c\colon \spc R \longrightarrow \Z_{\geq 0}$ is constant if and only if 
$\Gamma_{\small\kapp}$ is obtained from $\Gamma_K$
via a sequence of operations of type (I) and (II)a.
Otherwise, $c(\Gamma_{\small\kapp})>c(\Gamma_{K})$.

\end{prop}
\begin{proof}
Clearly, the weights on the graphs do not interfere with the complexity, as well as adding loops.

We now prove that, applying operation (II) the complexity remains the same only in case (II)a, and it increases strictly in case (II)b.
Call $\Gamma'$ the graph obtained by applying operation (II) to a vertex $v$ of degree $d$, and, accordingly to the notation of Proposition \ref{degenerations}, call 
$d_a$ and $d_b$ the vertex degrees of $v_a$ and $v_b$, respectively,  in $\Gamma'$. 
Clearly 
$$
d=d_a+d_b -2,
$$
Since we are adding a vertex and a component, the first Betti number of the graph, $b_1$, remains the same. 
However, the complexity can increase, but not decrease. 
Indeed, as observed above,  we have that $\Gamma=\Gamma'\cdot l$. 
So, by formula (\ref{fond}),
$$
c(\Gamma')=c(\Gamma)+c(\Gamma'-l)\geq c(\Gamma).
$$
Moreover, we see from the above formula that the complexity remains the same if and only if $c(\Gamma'-l)=0$, 
i.e., if and only if  $l$ 
is a disconnecting edge for $\Gamma'$. 
\end{proof}

\subsubsection*{The switching of two edges}
As we will often make use of another operation on graphs, it is convenient to describe it separately.
Let us consider a graph $\Gamma$, and fix two distinct edges $l$ and $m$, with associated vertices $v$, $v'$ and $w$, $w'$ respectively. Let us suppose that $v\not = w$ (but we do not ask that $v\not = v'$ or $w\not = w'$, i.e. $l$ and $m$ can be loops). 
We shall construct a new graph  $\Gamma'$ which has the same vertices as $\Gamma$, 
such that $\Gamma\setminus\{l,m\}$ is a subgraph of $\Gamma'$, and the edges $l$, $m$, are substituted by 
$l'$ and $m'$ as illustrated in figure below.

\bigskip

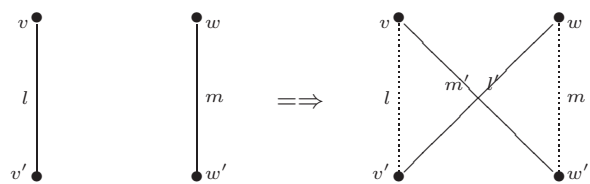
\centerline{
\begin{tabular}{rcl}
\xymatrix@R=2pc{
*{\bullet}  &   &*{\bullet} \\
& &\\
 *{\bullet}  \ar @{-}[uu]^{l}^>{v}^<{v'}   &   & *{\bullet} \ar @{-}[uu]_{m}_>{w}_<{w'}\\}
&
\xymatrix@R=2pc{
\\
=\Rightarrow\\
\\ 
\\}
&
\xymatrix@R=2pc{
*{\bullet} \ar @{-}[ddrr]^{l^\prime } &   & *{\bullet} \\
& &\\
 *{\bullet}  \ar @{.}[uu]^{l}^>{v}^<{v'}   \ar @{-}[uurr]^{m^\prime} &   & *{\bullet} \ar @{.}[uu]_{m}_>{w}_<{w'}\\}
\end{tabular}}

\bigskip

\noindent So, the new edges $l'$ and $m'$ connect $v$ and $w'$, $v'$ and $w$ respectively. 
We will call this process {\em the switching of  $l$ and $m$ with respect to $v$ and $w$\footnote{Note that this operation
 is the $\widetilde X$-transformation described in \cite{tsukui}.}.}
Note that the vertex degrees of $\Gamma$ and $\Gamma'$ are the same.

In general, this operation doesn't increase the complexity; counterexamples are easy to construct.
However, we shall prove that, for certain graph configurations, it does.

\begin{prop}\label{1}
Let $C$ be as in Theorem \ref{maxconf}. 
Then $C$ has no disconnecting nodes.
\end{prop}
\begin{proof}
Suppose that $C$ has a disconnecting node, which corresponds in $\Gamma=\Gamma_C$ to a disconnecting edge $r$.
Consider $l$, $m$ two  edges adjacent to $r$. 
The only case in which we cannot find two different edges, is if one of the vertex joined by  $r$, say $v$, has degree one.
In this case, as the curve $C$ is stable, necessarily $ g(v)\geq 1$; we therefore modify the graph by adding a loop on $v$, and 
decreasing by one the weight on $v$.
We still have a graph associated to a stable curve of genus $g$ with the same complexity, which we will call again $\Gamma$, and we can choose $l$ and $m$ as above.
Now, we perform the switching of  $l$ and $m$ w.r.t. $v$ and $w$:\footnote{this operation is the 
``slide transformation of $l$ and $m$ at $v$ and $w$ along $r$'' with the terminology of \cite{tsukui}.}

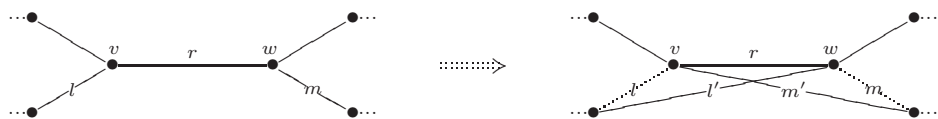
\bigskip

\centerline{
\xymatrix@R=1pc{
\ar @{}[r]|->{\ldots} &*{\bullet}  \ar @{-}[dr]^>{v} &  & &  & *{\bullet} \ar@{-}[dl]_>{w}& \ar @{}[l]|->{\ldots} &\ar @{}[r]|->{\ldots} &*{\bullet}  \ar @{-}[dr]^>{v} &  & &  & *{\bullet} \ar@{-}[dl]_>{w}& \ar @{}[l]|->{\ldots} \\
&& *{\bullet} \ar@{-}[rr]^{r} & & *{\bullet} & &\ar@{:>}[r] &&& *{\bullet} \ar@{-}[rr]^{r} & & *{\bullet} & \\
*{} \ar @{}[r]_>{}|->{\ldots}& *{\bullet}  \ar @{-}[ur]|-{l} &  & &  & *{\bullet} \ar @{-}[ul]|-{m}  & \ar @{}[l]^>{}|->{\ldots}&
*{} \ar @{}[r]_>{}|->{\ldots}& *{\bullet}  \ar @{.}[ur]|-{l} \ar @{-}[urrr]|-{l^\prime} &  & &  & *{\bullet} \ar @{.}[ul]|-{m} \ar @{-}[ulll]|-{m^\prime} & \ar @{}[l]^>{}|->{\ldots}\\}
}

\bigskip

Note that $\Gamma\cdot r=\Gamma' \cdot r$.
Now, $\Gamma'$ is the dual graph of a stable curve (indeed, we can think of the modification made as if we have 
resolved two nodes of the curve, and attached the components where they belong in a different way).

Now, if the edge $r$ is still a disconnecting edge for $\Gamma'$, then both $l$ and $m$ should be also disconnecting edges in $\Gamma$. In this case, we must have $g(v)\geq 1$ and we perform operation (I) in $v$. Then, if we choose $l$ to be the new loop, after the switching of $l$ and $m$ with respect to $v$, the edge $r$ will not be diconnecting anymore.

So, we can suppose that $r$ is no longer a disconnecting edge for $\Gamma'$, and it is immediate to check that we have 
introduced no new disconnecting edges, i.e. if $f$ is a disconnecting edge for $\Gamma'$, then it is also a disconnecting edge for $\Gamma$.

Now, applying Proposition \ref{comb} and observing  that $c(\Gamma-r)=0$, while $c(\Gamma'-r)>0$, we get:
$$
c(\Gamma')=c(\Gamma'\cdot r)+c(\Gamma'-r)=c(\Gamma\cdot r)+c(\Gamma'-r)>c(\Gamma\cdot r)=c(\Gamma).
$$
This proves the statement.
\end{proof}

\begin{rem}
\upshape{In \cite{DJ-IR}, sec.3, it is observed that a general graph (not necessarily regular) has to be
free of disconnecting edges in order to have maximal complexity. 
Proposition \ref{1}, as stated, could be indeed proven in the same way.
However, the operation used by the two authors is the switching 
of {\em one} edge, which modifies the vertex degrees; as we will see, it is necessary for us to remain in the class 
of cubic graphs, that's why we apply our switching.}
\end{rem}


\begin{prop}\label{2}
Let $C$ be as in Theorem \ref{maxconf}. 
Then $C$ is a loopless graph curve.
\end{prop}
\begin{proof}
From Proposition \ref{1} we can suppose that $C$ has no disconnecting nodes.
Let $\Gamma$ be the weighted graph associated to $C$ and suppose $\Gamma$ is not a loopless trivalent graph.
Using the degeneration operations described in Definition \ref{degenerations},
we will describe an algorithm in two steps that, given $\Gamma$ with no disconnecting nodes,
produces a loopless graph $\Gamma_3$ 
with strictly bigger complexity than $\Gamma$, such that $b_1(\Gamma_3)=g$ and each vertex has degree $3$ and weight  $0$. 
Therefore $\Gamma_3$ is the dual graph associated to a loopless graph curve of genus $g$.

\medskip

\noindent FIRST STEP (Reduction to a curve having only rational components): 
We replace any vertex $v_i$ with strictly positive weight $g(v_i)$ with a bouquet given by a vertex of weight $0$ and $g(v_{i})$ 
loops attached to it. 
This is a reiterate application of operation (I) and does not change the arithmetic genus, nor the complexity. 
Call $\Gamma_1$ the graph obtained in this way.

\medskip

\noindent SECOND STEP (Reduction to a loopless graph curve):
suppose that  in $\Gamma_1$ there is a vertex $v$ with degree $d$  greater or equal to $4$. 
As $\Gamma_1$ is obtained from $\Gamma$ by adding loops, it has no disconnecting edges. 
We can therefore perform operation (II)b on $v$, with the request that $\deg(v_a)\geq 3$, and $\deg(v_b)\geq 3$
(i.e. remaining in the class of graphs of stable curves).
As proved in Proposition \ref{increasing}, this operation increases strictly the complexity.
Moreover, the degrees of $v_a$ and $v_b$ are strictly smaller than the degree of $v$.
We perform this transformation until all the vertices of the new graph have exactly degree 3. 
Call $\Gamma_2$ the resulting graph. 
So by the genus formula $\Gamma_2$ has $2g-2$ vertices of weight 0, $3g-3$ edges, and $b_1(\Gamma_2)=g$; 
equivalently, $\Gamma_2$ is the graph associated to a graph curve.
Moreover, $\Gamma_2$ has no disconnecting edges,  and no loops. 
Indeed, if there were a loop, $\Gamma_2$  would contain a subgraph  of the form:

\bigskip

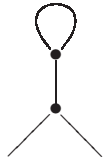
\begin{equation*}
\xymatrix@=1.2pc{
&*{\bullet} \ar@{-}@(ur,ul) \ar@{-}[d]\\
&*{\bullet} & \\
\ar@{-}[ur] && \ar@{-}[ul] & .
}
\end{equation*}
hence with a disconnecting node.
Moreover, by what observed above, $c(\Gamma_2)> c(\Gamma_1)$.
\end{proof}



\begin{ex}
\upshape{
Let us consider the dual  graph of a smooth genus $3$ curve; this is just one vertex with weight $3$.
Let us apply the above algorithm to this graph.
With the first step, we obtain  the following \av bouquet" of $3$ loops,

\medskip

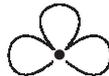
\begin{equation*}
\xymatrix@R=1pc{
*{\bullet} \ar@{-}@(ul,dl) \ar@{-}@(ur,dr) \ar@{-}@(ur,ul)
}
\end{equation*}

\medskip

which is the graph of a stable curve with a single irreducible rational component with $3$ nodes. 
As the irreducible component  has geometric genus $0$, we apply the second step. 
We have several choices. 
So, we can obtain either one of the following configurations.

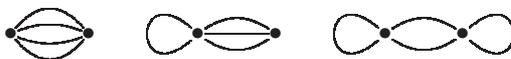
\begin{equation*}
\begin{tabular}{cccccccccc}
\xymatrix{
*{\bullet}  \ar @{-} @/_.3pc/[r] \ar @{-}@/_.8pc/[r]/& *{\bullet} \ar@{-}@/_.3pc/[l] \ar@{-}@/_.8pc/[l]
}
&{}&{}&
\xymatrix{
*{\bullet} \ar@{-}@(ul,dl) \ar @{-} @/_/[r] \ar@{-}[r] & *{\bullet} \ar@{-}@/_/[l]
}
&{}&{}&
\xymatrix{
*{\bullet} \ar @{-}@(ul,dl) \ar @{-}@/_/[r] & *{\bullet} \ar @{-}@/_/[l] \ar @{-}@(ur,dr)
}
\end{tabular}
\end{equation*}
}

Notice that all these graphs have still vertices with degree greater then 3. So, we therefore must go on applying step 2. 
For example, one of the possibilities for the first graph would be the following.

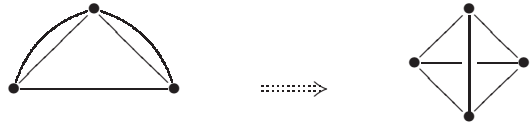
\begin{equation*}
\begin{tabular}{lcr}
\xymatrix{
& *{\bullet} \ar@{-}[dl] \ar@{-}@/_/[dl] \ar@{-}@/^/[dr] \ar@{-}[dr]& \\
*{\bullet} \ar@{-}[rr] & & *{\bullet} & \ar@{:>}[r] &
}
&
&\xymatrix@=1.2pc{
&*{\bullet}\ar@{-}[dd] \ar@{-}[dl] \ar@{-}[dr] &\\
*{\bullet} \ar@{-}[rr] |!{[r]}\hole & &*{\bullet}\\
&*{\bullet}\ar@{-}[ru]\ar@{-}[lu] & 
}
\end{tabular}
\end{equation*}

If we start from the third graph, we can perform the following chain of transformations
(indeed, in this case, these operations are forced).

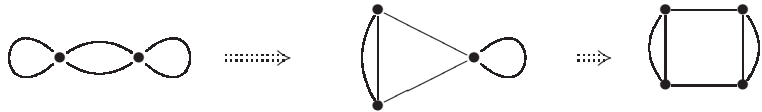
\begin{equation*}
\begin{tabular}{lcccr}
\xymatrix@R=1pc{
{}\\
 *{\bullet} \ar @{-}@(ul,dl) \ar @{-}@/_/[r] & *{\bullet} \ar @{-}@/_/[l] \ar @{-}@(ur,dr)&\ar@{:>}[r]&
}
&&
\xymatrix@=1pc{
*{\bullet} \ar@{-}[dd]\ar@{-}@/_.5pc/[dd]\ar@{-}[drr]&&&&\\
&&*{\bullet} \ar@{-}@(ur,dr)&&\ar@{:>}[r]&\\
*{\bullet} \ar@{-}[urr]&& 
}
&
\xymatrix{
*{\bullet} \ar@{-}[r] \ar@{-}[d] \ar@{-}@/_/[d] & *{\bullet} \ar@{-}[d] \ar@{-}@/^/[d]\\
*{\bullet} \ar@{-}[r] & *{\bullet}
}
\end{tabular}
\end{equation*}
\end{ex}

As also the simple example above shows, the algorithm of Proposition \ref{2} doesn't give a unique output.
On the contrary, if we start from a bouquet of $g$ loops we can obtain \emph{all} the possible cubic graphs without disconnecting edges  using step 2,  as one can see in the following way. 
Consider a cubic graph without disconnecting edges  $\Sigma$ with $b_1(\Sigma)=g$. 
Let $T$ be a  spanning tree for $\Sigma$. To contract all the edges of $T$ is the reverse operation to step 2 of the algorithm, and the result is the bouquet of $g$ loops.

\begin{proof} ({\it of Theorem \ref{maxconf}})
By Proposition \ref{1} and \ref{2} we can suppose that $\Gamma=\Gamma_C$  is trivalent,
loopless and without disconnecting edges.
Hence, the only possible multiple edges are of the form

\medskip

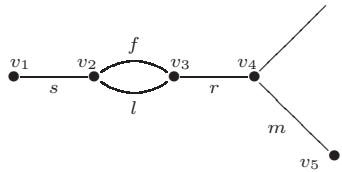
\centerline{
\xymatrix{ &&&&\ar @{-}[dl]\\
*{\bullet}  \ar @{-}[r]_s^>{v_2}^<{v_1}&*{\bullet}  \ar @{-}@/_/[r]_l& *{\bullet} \ar @{-}@/_/[l]_f & *{\bullet}  \ar @{-}[l]^r_>{v_3}_<{v_4}&\\
&&&& *{\bullet}\ar @{-}[ul]^m^<{v_5}\\
}
}

\bigskip

Note that $l$ and $m$ have to be distinct as well as $v_1$ and $v_4$ (otherwise $\Gamma$ would have disconnecting edges).
With the notations adopted in the figure, we apply the switching of  $l$ and $m$ w.r.t. $v_3$ and $v_4$:

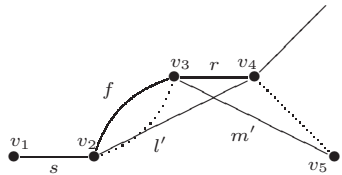
\centerline{
\xymatrix{ &&&&\ar @{-}[dl]\\
& & *{\bullet} \ar @{-}@/_/[dl]_{f} & *{\bullet}  \ar @{-}[l]_r_>{v_3}_<{v_4}&\\
*{\bullet}  \ar @{-}[r]_s^<{v_1}^>{v_2}&*{\bullet}  \ar @{.}@/_/[ur]_{l^\prime}\ar @{-}[urr]&&&*{\bullet} \ar @{.}[ul]\ar @{-}[ull]^{m^\prime}^<{v_5}\\
}
}

\bigskip

\bigskip

We have constructed a new trivalent graph $\Gamma'$ which has one less couple of multiple edges: indeed, 
it is immediate to check that we have not created any other multiple edges.

Let us prove that $c(\Gamma ')>c(\Gamma)$.
We prove it by giving an injective map from the spanning trees of $\Gamma$ to the spanning trees of $\Gamma'$, as follows.

Let $T$ be a spanning tree of $\Gamma$. If $r\in T$, the switching of $\Gamma$ we performed transforms $T$ into a spanning tree $T^\prime$ of $\Gamma^\prime$, with $r\in \Gamma^\prime$.

Now, suppose $r \not\in T$. We shall distinguish between $4$ different situations.
\begin{itemize}
\item $l,m\not\in T$ ($\Rightarrow f\in T$). $T$ is transformed into a tree $T^\prime$ of $\Gamma^\prime$ with $l^\prime\not\in T^\prime,m^\prime\not\in T^\prime,r\not\in T^\prime$ and $f\in T^\prime$. 
\item $l\in T$, $m\not\in T$ ($\Rightarrow f\not\in T$). Then, in $T$, the path connecting $v_2$ and $v_4$ passes by $l$. So, the switching of $l$ and $m$ along $r$ restricted to $T$ is not a spanning tree of $\Gamma^\prime$ since $l^\prime$ would create a cycle containing $l^\prime$ and $v_3$ and $v_5$ get disconnected. So, we consider the subgraph of $\Gamma^\prime$ obtained by the union of $m^\prime$ with the switching of $T$ minus $l^\prime$. It is easy to see that this is a spanning tree $T^\prime$ of $\Gamma^\prime$, with $l^\prime\not\in T^\prime$, $m^\prime\in T$, $r\not\in T^\prime$ and $f\not\in T^\prime$.
\item $l\not\in T$, $m\in T$ ($\Rightarrow f\in T$). In this case we must do a further distinction. 
\begin{itemize}
\item In $T$ the path connecting $v_2$ and $v_4$ does not contain $m$. 

In this case $T$ is transformed into a spanning tree $T^\prime$ of $\Gamma^\prime$ with $l^\prime\not\in T^\prime$, $m^\prime\in T^\prime$, $r\not\in T^\prime$ and $f\in T^\prime$.
\item In $T$ the path connecting $v_2$ and $v_4$ contains $m$.

In this case the switching we performed restricted to $T$ is not a spanning tree of $\Gamma^\prime$, since $m^\prime$ creates a cycle and $v_4$ and $v_2$ get disconnected. So, we consider $T^\prime$ as the union of $l^\prime$ with the switching of $\Gamma$ restricted to $T$ minus $m^\prime$. Then $T^\prime$ is a spanning tree of $\Gamma^\prime$ with $l^\prime\in T^\prime$, $m^\prime\not\in T^\prime$, $r\not\in T^\prime$ and $f\in T^\prime$.
\end{itemize}
\item $l,m\in T$ ($\Rightarrow f\not\in T$). Also in this case we have to distinguish between the following situations:
\begin{itemize}
\item In $T$ the path connecting $v_2$ and $v_4$ does not contain $m$.

Then the switching of $l$ and $m$ along $r$ restricted to $T$ is not a spanning tree of $\Gamma^\prime$ since it creates a cycle containing $l^\prime$ and $v_4$ and $v_5$ get disconnected. So if we consider $T^\prime$ as the union of $f$ with the switching of $T$ minus $s$, then $T^\prime$ is a spanning tree of $\Gamma^\prime$ with $l^\prime\in T^\prime$, $m^\prime\in T^\prime$, $r\not\in T^\prime$ and $f\in T^\prime$.
\item In $T$ the path connecting $v_2$ and $v_4$ contains $m$.

In this case the switching we performed restricted to $T$ is a spanning tree $T^\prime$ of $\Gamma^\prime$ with $l^\prime\in T^\prime$, $m^\prime\in T^\prime$, $r\not\in T^\prime$ and $f\not\in T^\prime$.
\end{itemize}
\end{itemize}

So, the association of each spanning tree $T$ of $\Gamma$ with $T^\prime$ gives an injective map from the spanning trees of $\Gamma$ into the spanning trees of $\Gamma^\prime$.

Now, to prove that, indeed strict inequality holds, we observe that $\Gamma^\prime$ has at least one spanning tree $T^\prime$ such that $l'$, $f$ and $r$ are not in $T^\prime$ ($\Rightarrow m^\prime\in T^\prime$), and this kind of spanning trees are not in the image of the map constructed above. 
\end{proof}




\subsection{Combinatorial version}\label{combi}
We can give a purely combinatorial translation of the above results, as follows.
Recall that a graph $\Gamma$ is said to be {$k$-connected} if  for any  set $S$ of $k-1$ edges of 
 $\Gamma$, $\Gamma\setminus S$ is still connected. 
 $\Gamma$ is said to be {\it strictly} $k$-connected if it is $k$-connected but not $k+1$-connected. 
So, a graph is $1$-connected  if and only if it is connected; it is strictly $1$-connected if and only if it is connected and 
 it has at least one disconnecting node.
 Of course, a trivalent graph can never be $4$-connected (any triple of edges incident in one vertex is a disconnecting set).
 
 Given a positive integer $g$, let $\mathcal G_g$ be the set of all weighted multigraphs of genus $g$ satisfying 
 condition (\ref{stability}). Theorem \ref{maxconf} can now be written in the following way.
 
 \begin{teo}\label{combinmaxconf}
 The graphs $\Gamma\in \mathcal G_g$ reaching the maximal complexity
 are trivalent, simple and $2$-connected.
 \end{teo}

We shall call a simple trivalent  graph a {\em cubic graph}.  


\section{Bounds on the maximal complexity of stable curves}\label{ub}

The question of an upper bound on the complexity depending on the genus has several geometrical meanings: for instance, it implies that there is a bound for the irreducible components of the compactified Jacobian of a stable curve, and a bound for the group of components of the N\'eron model of the relative Jacobian of a family of curves having $C$ as special fiber; see also Remark \ref{confrontoL} .
Let us define 
$$\psi(g):=\mbox{max}\{c(C), C \mbox{ stable curve of genus } g\}.$$

\begin{rem}\upshape{
Notice that $\psi$ is strictly increasing with respect to the genus $g$. 
In fact, if $C$ is a curve of genus $g$ such that $\psi(g)=|\Delta_C|$, let $C^\prime$ be the stable curve obtained from $C$ by adding an extra node connecting two of the components of $C$. Then, using the above formula for the genus of $C^\prime$, one gets $g_{C^\prime}=g+1$ and, clearly, $c(C^\prime)>c(C)$.}
\end{rem}

According to the results of the preceding section, we know that curves achieving maximal complexity are in particular graph curves. 
Using this fact, we can find a first rough bound:
$$\psi(g)\le \binom{3g-3}{g}\le 2^{3g-4}.$$  
For the proof, see \cite{lorjac}, Lemma 2.7.
Note that  the first inequality follows immediately from the fact that any spanning tree in a cubic
graph of order $2g-2$ is obtained by removing $g$ edges.

Applying the results on the surjectivity of the ``Abel-Jacobi map'' in \cite{B-N},
we could derive a better bound: $$\psi(g)\le \binom{2g-2}{g}.$$

Both these bounds are not sharp, even for low genus.
Using a result of Biggs \cite{biggs}, we can prove the following.
\begin{equation}\label{disbiggs}
c(\Gamma)\leq  \frac{1}{2g-2}\left(\frac{6g-6}{2g-3}\right)^{2g-3}=:\alpha(g).
\end{equation}

For $g=3$, this bound is optimal.
However, the bound is not asymptotically sharp.

The complexity of $k$-regular graphs has been studied also  by McKay \cite{McK}, and Chung-Yau \cite{C-Y}.
We can apply their results obtaining a sharper bound. 
\begin{teo}\label{bound}
Let $C$ be a stable curve of genus g, and let $\Gamma$ be its dual graph. 
Then 
\begin{equation}\label{disCY}
c(\Gamma)\leq  \frac{2 \ln{(2g-2)}}{3(2g-2) \ln{(9/8)}}
\exp{\left(\frac{12}{\sqrt\pi}\left(\frac{\ln{9/8}}{\ln{(2g-2)}}\right)^{5/2}\right)}
\left(\frac{4}{\sqrt{3}}\right)^{2g-2}=:\beta(g).
\end{equation}
Moreover, this bound is asymptotically sharp for $g\gg 0$.
\end{teo}
\begin{proof}
By making explicit the computations in Theorem 4 of \cite{C-Y}, we obtain the bound (\ref{disCY})
for $3$-regular graphs with $2g-2$ vertices. The thesis is now straightforward because of the reduction to
graph curves of Theorem \ref{maxconf}.
McKay proves in \cite{McK} that there is a sequence of cubic graphs 
$\{\Gamma_i\}_{i\in\mathbb N}$ with increasing orders $n_i$'s, such that, if we set 
\begin{equation}\label{tau}
\tau(\Gamma_i):=\left(n_ic(\Gamma_i)\right)^{\frac{1}{n_i}},
\end{equation}
then $\tau(\Gamma_i)\longrightarrow \frac{4}{\sqrt 3}$ for $i\rightarrow \infty$.
Indeed, he proves much more than the existence, but that a ``random'' sequence of cubic
graphs satisfies this property, see also \cite{lyons}.
Hence, the constant $ \frac{4}{\sqrt 3}$ in the bound (\ref{disCY}) is the best possible.
\end{proof}

In the next picture we describe the graphs with maximal complexity for $g=4$, $g=5$ and $g=6$, respectively.

\bigskip

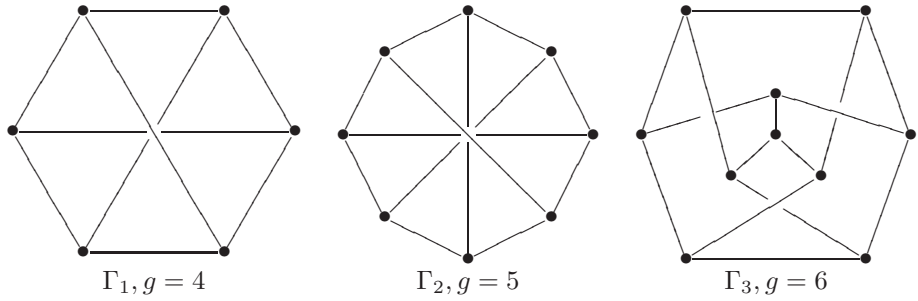
\centerline{
\begin{tabular}{ccc}
\xymatrix@C=1.7pc@R=3.3pc{
& *{\bullet} \ar@{-}[rr] \ar@{-}[ddrr] & & *{\bullet} \ar@{-}[ddll]|!{[ll];[dd]}\hole &\\
*{\bullet} \ar@{-}[rrrr]|!{[rr]}\hole \ar@{-}[ur] \ar@{-}[dr]  & & & & *{\bullet} \ar@{-}[ul] \ar@{-}[dl] \\
& *{\bullet} \ar@{-}[rr] & & *{\bullet}
}
&
\xymatrix@=.8pc{
& & & *{\bullet} \ar@{-}[dll] \ar@{-}[drr] \ar@{-}[dddddd]|!{[ddd]}\hole & & & \\
& *{\bullet} \ar@{-}[ddl] \ar@{-}[ddddrrrr] & & & & *{\bullet} \ar@{-}[ddr]
\ar@{-}[ddddllll]|!{[llll];[dddd]}\hole & \\
&&&&&&\\
*{\bullet} \ar@{-}[rrrrrr]|!{[rrr]}\hole & & &&& & *{\bullet} \\
&&&&&&\\
& *{\bullet} \ar@{-}[uul] & & & & *{\bullet} \ar@{-}[uur] & \\
& & &  *{\bullet} \ar@{-}[ull] \ar@{-}[urr]
}
&
\xymatrix@C=.9pc@R=.8pc{
&  *{\bullet} \ar@{-}[rrrr] \ar@{-}[dddl] \ar@{-}[ddddr] &&&& *{\bullet} \ar@{-}[dddr] \ar@{-}[ddddl]|!{[ddd];[lll]}\hole & \\
&&&&&&\\
&  & & *{\bullet} \ar@{-}[d] \ar@{-}[dlll]|!{[d];[lll]}\hole \ar@{-}[drrr]&&&\\
*{\bullet} \ar@{-}[dddr] &&  &*{\bullet} \ar@{-}[dl] \ar@{-}[dr]& && *{\bullet} \ar@{-}[dddl] \\
&&*{\bullet} \ar@{-}[ddrrr]|!{[d];[rrr]}\hole &&*{\bullet} \ar@{-}[ddlll]&&\\
&&&&&&\\
& *{\bullet} \ar@{-}[rrrr] &&&&*{\bullet}
}\\
$\Gamma_1 , g=4$
&
$\Gamma_2 , g=5$
&
$\Gamma_3 , g=6$
\end{tabular}
}

\bigskip

Note that $\Gamma_3$ is the famous Peterssen graph. 
These graphs are proved to be of maximal complexity among simple cubic graphs of their respective orders in \cite{DJ-IR}, sec.5. 



\begin{rem}\label{confrontoL}
\upshape{As already noticed, the geometric meaning of this result is that it gives a bound on the group of connected components of the N\'eron model of the degree-$d$ Picard variety for families of {\em stable}
curves (\cite{BLR}, Theorem 1, sec.9.6), and as well on the number of irreducible components of the fibres of the scheme
$P_g^d$ constructed in \cite{capneron} and of $\overline{P_{d,g}}$ of \cite{cap}. 
Of course this is also a bound on the cardinality of the group of components of the N\'eron models of the Jacobians of such families.
However, notice that this bound is different from the ones found  by Lorenzini in \cite{lorner} (see also \cite{BLR}, 
Theorem 9 sec.6.9). 
Indeed, given any strictly henselian discrete valuation ring $R$, with algebraically closed residue field $\kapp$ and field of fractions $K$,  any regular family of curves $f\colon X\longrightarrow \mbox{spec }R$ such that the {\em general} fiber $X_K$ is of genus $g$, and the Jacobian $J_K$  has {\em potential good reduction}, Lorenzini finds a bound, depending only on 
$g$, for the group of components of the N\'eron model of $f$.
We obtain instead a bound for the group of components of the N\'eron model of \emph{any} smooth family $f\colon X\rightarrow \mbox{spec }R$ such that the \emph{closed} fiber $X_{\kappy}$ is a  stable curve.
So one could say that his approach is dynamic, in the sense that it depends on the generic fiber, while our approach is static, as it starts from a given curve.
Moreover, the boundedness comes on the one hand from the assumption of potential good reduction on the general fiber, on the other from the assumption of stability of the special one.
As it is proven in \cite{lorgraphs}, a  family has potential good reduction if and only if the closed fiber has a {\em tree} as dual graph, 
hence complexity 1.
}
\end{rem}


\subsection{Example and lower bound}\label{lower}




We present here an example of a family of cubic graphs of increasing order and increasing complexity.
In particular, this gives explicit lower bounds on $\psi(g)$ (even though not sharp, given the above result).

Let $\Gamma_m$ be the graph with $4m$ vertices of valency $3$ formed from  $m$ 
pairwise disjoint graphs $G_i$ of the following form:

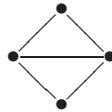
\centerline{
\begin{tabular}{c}
\xymatrix@=1pc{
& *{\bullet} \ar@{-}[dl] \ar@{-}[dr]& \\
*{\bullet} \ar@{-}[rr] && *{\bullet} \\
& *{\bullet} \ar@{-}[ul] \ar@{-}[ur]
}
\end{tabular}
}

\noindent by adding $m$ edges $l_1,\ldots ,l_m$ to link them as a ring, as shown in the figure below for $m=6$. 
Clearly, $\Gamma_m$ is the dual graph of a graph curve of genus $g=2m+1$.

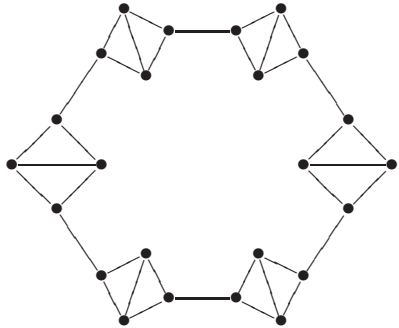
\centerline{
\xymatrix@=.2pc{
&&&&&*{\bullet} \ar@{-}[drr] \ar@{-}[ddl] \ar@{-}[dddr] &&&&&&& *{\bullet} \ar@{-}[dll] \ar@{-}[ddr] \ar@{-}[dddl]&&& \\
&&&&&&&*{\bullet} \ar@{-}[ddl] \ar@{-}[rrr] &&&  *{\bullet} \ar@{-}[ddr] &&&&& \\
&&&&*{\bullet} \ar@{-}[drr] \ar@{-}[dddll] &&&&&&&&& *{\bullet} \ar@{-}[dll] \ar@{-}[dddrr]\\
&&&&&&*{\bullet} &&&&& *{\bullet} &&&& \\
\\
&&*{\bullet} \ar@{-}[ddll] \ar@{-}[ddrr] &&&&&&&&&&&&& *{\bullet} \ar@{-}[ddrr]\ar@{-}[ddll]& \\
\\
*{\bullet} \ar@{-}[rrrr] \ar@{-}[ddrr] &&&& *{\bullet} \ar@{-}[ddll] &&&&&&&&& *{\bullet}  \ar@{-}[ddrr] \ar@{-}[rrrr] &&&& *{\bullet} \ar@{-}[ddll]\\
\\
&&*{\bullet} \ar@{-}[dddrr] &&&&&&&&&&&&& *{\bullet} \ar@{-}[dddll] & \\
\\
&&&&&& *{\bullet} \ar@{-}[dll] \ar@{-}[ddr] \ar@{-}[dddl] &&&&& *{\bullet} \ar@{-}[drr] \ar@{-}[ddl] \ar@{-}[dddr]\\
&&&& *{\bullet} \ar@{-}[ddr] &&&&&&&&& *{\bullet} \ar@{-}[ddl]\\
&&&&&&&*{\bullet} \ar@{-}[dll] &&& *{\bullet} \ar@{-}[lll] \ar@{-}[drr]\\
&&&&& *{\bullet} &&&&&&& *{\bullet}
}
}

\vspace{.1cm}


\begin{prop}
$
c(\Gamma_m)=2m8^m.
$
\end{prop}
\noindent {\bf Proof}: To prove this formula we will proceed by hands counting the number of spanning trees for $\Gamma_m$.
Notice that the complexity of the subgraphs $G_i$'s is $8$. 

Choosing one of the edges $l_i$, the number of spanning trees not containing it is $8^m$. 
Indeed, all the others $l_j$'s have to be included in any spanning tree, while for any $G_k$, we have to count the $8$ possibilities for the spanning trees. So, we have $m8^m$ of these.

On the other hand, if $T$ is a spanning tree that contains all the $l_i$, then there is one and only one $j$ such that $T\cap G_j$ is disconnected, and it has to be of  one of the $8$ following forms:

\vspace{.2cm}
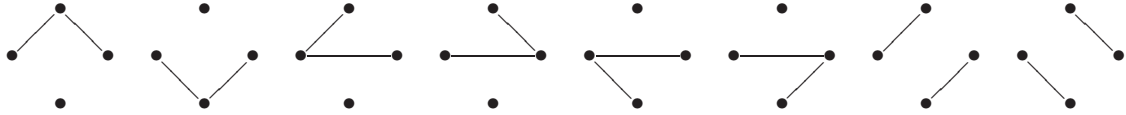
\begin{tabular}{cccccccc}
\xymatrix@=1pc{
& *{\bullet} \ar@{-}[dl] \ar@{-}[dr]& \\
*{\bullet}  && *{\bullet} \\
& *{\bullet}
}
&
\xymatrix@=1pc{
& *{\bullet}& \\
*{\bullet} && *{\bullet} \\
& *{\bullet} \ar@{-}[ul] \ar@{-}[ur]
}
&
\xymatrix@=1pc{
& *{\bullet} \ar@{-}[dl] & \\
*{\bullet} \ar@{-}[rr] && *{\bullet} \\
& *{\bullet} 
}
&
\xymatrix@=1pc{
& *{\bullet}  \ar@{-}[dr]& \\
*{\bullet} \ar@{-}[rr] && *{\bullet} \\
& *{\bullet} 
}
&
\xymatrix@=1pc{
& *{\bullet} & \\
*{\bullet}  \ar@{-}[rr]&& *{\bullet} \\
& *{\bullet} \ar@{-}[ul] 
}
&
\xymatrix@=1pc{
& *{\bullet} & \\
*{\bullet} \ar@{-}[rr] && *{\bullet} \\
& *{\bullet} \ar@{-}[ur]
}
&
\xymatrix@=1pc{
& *{\bullet} \ar@{-}[dl] & \\
*{\bullet}  && *{\bullet} \\
& *{\bullet} \ar@{-}[ur]
}
&
\xymatrix@=1pc{
& *{\bullet}  \ar@{-}[dr]& \\
*{\bullet}  && *{\bullet} \\
& *{\bullet} \ar@{-}[ul] 
}
\end{tabular}

\vspace{.2cm}

For $i\not =j$, $T\cap G_i$ is a spanning tree for $G_i$ as above.
Therefore, we have $8m8^{m-1}=m8^m$ possible spanning trees of this form.

So, summing all up, we have $2m8^m$ spanning trees, as required.\hfill $\Box$

\bigskip

Notice that $$c(\Gamma_m) =(g-1)8^{\frac{g-1}{2}}=(g-1)(2\sqrt{2})^{g-1}.$$
So, for odd $g$, we have this lower bound on $\psi$. 
For even $g$, remembering that $\psi$ is an increasing function of the genus, we have
$\psi(g)>\psi(g-1)$. 
Hence, $\psi$ is bounded from below by $(g-2)(2\sqrt{2})^{g-2}$.








\section{Further results and conjectures}\label{congetture}

\subsubsection*{3-connectedeness}

Note that a trivalent graph with loops, or with double edges is strictly $2$-connected. 
With Theorem \ref{MAX}, we have therefore excluded the strictly $1$-connected case, 
and some of the cases of strict $2$-connectedness. 
>From this, and from the observation of graphs with maximal complexity for low genus, 
it is natural to ask

\begin{conj}\label{connectivity}
A cubic graph with maximal complexity is $3$-connected.
\end{conj}
This seems to be generally believed, also for bigger classes of graphs (see e.g.\cite{DJ-IR} sec.3), 
but no proof is known.
Let us note that $3$-connected graph curves have also interesting geometrical properties; 
indeed, they are precisely those graph curves for which the canonical bundle is very ample, yielding an embedding morphism \cite{B-E}.

With similar techniques to the ones used in Section \ref{MAX}, we can prove a partial result.
\begin{prop}\label{bah}
Let $C$ be as in Theorem \ref{maxconf}. 
Then 
\begin{enumerate}
\item $\Gamma_C$ has no couple of disconnecting edges that lay on a cycle of length $\leq 6$.
\item  $\Gamma_C$ has no couple of disconnecting edges such that at least one of them is adjacent to a cycle of length $\leq 4$.
\end{enumerate}
\end{prop}

\begin{cor}\label{g<10}
Conjecture \ref{connectivity} holds for $g\leq 8$.
\end{cor}

>From a result on the so-called ``Abel-Jacobi map'' for graphs proved on \cite{B-N} (Theorem 1.8),
we can derive the following 
\begin{prop}[Baker-Norine]
If $\Gamma$ is any $k$-connected graph of order $n$, 
$$
c(\Gamma)\geq \binom{n}{k-1}.
$$
\end{prop}
This proves that the complexity grows with the connectivity.
In the case of cubic graphs, however, the bound obtained for $3$-connected graphs is just 
quadratic in the order.

\subsubsection*{Maximal girth}
The explicit sequences of cubic graphs with large complexity 
have all {\em large girth}. 
Let us explain the terminology.
The {\em girth} of a graph is the length of the shortest cycle.
A simple argument shows that in a cubic graph the girth cannot exceed $2\ln n/\ln2$ ($n$
being the order of the graph).
Hence, the girth of a family of cubic graphs can grow at most as the rate of the logarithm of 
the number of vertices.

\begin{defi}\label{largegirth}
A sequence $\{\Gamma_i\}_{i\in \mathbb N}$  of cubic graphs with increasing 
orders $n_i$ is called {\em sequence of large girth} if
$$
\lim_{i\to \infty}\frac{\ln(n_i)}{\ln2 \,\, girth(\Gamma_i)} \,\,\,\mbox{ is finite. }
$$
\end{defi}
McKay proves in \cite{McK} that sequences with large girth satisfy  condition (\ref{tau}).

\begin{conj}\label{girth}
A sequence $\{\Gamma_i\}_{i\in \mathbb N}$  of cubic graphs with maximal complexity has large girth. 
\end{conj}
See \cite{DJ-IR}, sec.3 for an heuristic argument supporting this conjecture.
 

From a result of McKay we can derive the following property of sequences of trivalent graphs of maximal complexity.
Given a graph $\Gamma$, and a positive integer $m$, let $C_{\Gamma}(m)$ be the number of 
cycles  of length $\leq m$  in $\Gamma$.

\begin{teo}\label{cicli}
Let $m$ be a positive integer. 
For any $\epsilon>0$ there exist a positive integer $j$ such that: 
given any infinite sequence of trivalent graphs
 $\{\Gamma_i\}_{i\in \mathbb N}$  such that 
$\Gamma_i$ has $2i$ vertices, and reaches the maximal complexity between cubic graphs of order $2i$,
the following inequality holds
$$
\frac{C_i(m)}{2i}\leq \epsilon \,\,\mbox{ for any }i\geq j,
$$
where $C_i(m)=C_{\Gamma_i}(m)$.
\end{teo}
\begin{proof}
Theorem 4.5 of \cite{McK} states that, given a sequence of cubic graphs $\{\Gamma_i\}_{i\in \mathbb N}$ 
 satisfying (\ref{tau}), then, for each fixed $m$,
 $$
\frac{C_i(m)}{n_i}\rightarrow 0 \,\,\mbox{ for } i\to \infty.
$$
This does not directly imply the statement, which is uniform, in the sense that the constant $j$ depends only on $m$ and not 
on the chosen sequence.
However, let us argue as follows: given $m$, consider a sequence of cubic  graphs with maximal complexity {\em and} with
maximal cardinality of $C(m)$.
Clearly this sequence satisfies condition (\ref{tau}).
Now, by assumption, the statement holds uniformly  on every sequence of graphs of maximal complexity.
\end{proof}
Note that a sequence of graphs of large girth satisfies trivially the condition of Theorem \ref{cicli}.


\subsubsection*{Further speculations}

Other open questions are: 
\begin{itemize}
\item Are cubic graphs of maximal complexity 
{\em strongly regular}?
\item Do cubic graphs of maximal complexity have {\em maximal automorphism group}? 
\item Is  a cubic graph with maximal complexity and given genus {\em unique}?
\item Let $C$ be a graph curve  with maximal complexity; suppose that it is indeed $3$-connected. 
Then the canonical map is an embedding realizing $C$ as a configurations of lines in $\mathbb{P}^{g-1}$.
We conjecture that these lines are in ``general position'', in a sense that can be made precise.
This property seems to be connected with the large girth property.
\item   Let $\Gamma$ be a cubic graph with maximal complexity; suppose that it is indeed $3$-connected. 
We think that the ``Clifford index'' of this graph, as defined for instance by Bayer and Eisenbud, has to be maximal as well. 
\item Let $\Gamma$ be a cubic graph with maximal complexity; suppose that it is indeed $3$-connected. 
We wonder if  any cutset of $3$ edges is made of adjacent edges.
\end{itemize}

\bibliographystyle{amsalpha}
  \nocite{*}
  \bibliography{complexity}

\providecommand{\bysame}{\leavevmode\hbox to3em{\hrulefill}\thinspace}
\providecommand{\MR}{\relax\ifhmode\unskip\space\fi MR }
\providecommand{\MRhref}[2]{%
  \href{http://www.ams.org/mathscinet-getitem?mr=#1}{#2}
}
\providecommand{\href}[2]{#2}
\begin{thebibliography}{BdlHN97}

\bibitem[Art91]{ar}
M.~Artin, \emph{Algebra}, Prentice Hall, 1991.

\bibitem[Art04]{artamkincan}
I.~V. Artamkin, \emph{Canonical mappings of punctured curves with the simplest
  singularities}, Mat. Sb. \textbf{195} (2004), no.~5, 3--32.

\bibitem[Art05]{artamkingen}
\bysame, \emph{Generating functions of modular graphs, and the {B}urgers
  equation}, Mat. Sb. \textbf{196} (2005), no.~12, 3--32.

\bibitem[BdlHN97]{tatiana}
R.~Bacher, P.~de~la Harpe, and T.~Nagnibeda, \emph{The lattice of integral
  flows and the lattice of integral cuts on a finite graph.}, Bull. Soc. Math.
  France \textbf{125} (1997), no.~2, 167--198.

\bibitem[BE91]{B-E}
D.~Bayer and D.~Eisenbud, \emph{Graph curves.}, Adv. Math. \textbf{86} (1991),
  no.~1, 1--40, {W}ith an appendix by {S}ung {W}on {P}ark.

\bibitem[Ber70]{berge}
C.~Berge, \emph{Graphs and hypergraphs}, North-Holland, 1970.

\bibitem[Big74]{biggs}
N.~Biggs, \emph{Algebraic {G}raph {T}heory}, Cambridge University Press, 1974.

\bibitem[Big99]{biggschip}
\bysame, \emph{Chip-{F}iring and the {C}ritical {G}roup of a {G}raph}, Journal
  of Algebraic Combinatorics \textbf{9} (1999), 25--45.

\bibitem[BLR90]{BLR}
S.~Bosch, W.~L\"uktebohmert, and M.~Raynaud, \emph{N\'eron models}, Ergebnisse
  der Mathematik, no.~21, Springer-Verlag, 1990.

\bibitem[BMS06]{BMS}
S.~Busonero, M.~Melo, and L.~Stoppino, \emph{Combinatorial aspects of stable
  curves}, Le Matematiche \textbf{LXI} (2006), no.~I, 109--141.

\bibitem[BN07]{B-N}
M.~Baker and S.~Norine, \emph{Riemann-{R}och and {A}bel-{J}acobi theory on a
  finite graph}, Adv. Math. \textbf{215} (2007), no.~2, 766--788.

\bibitem[Bol98]{bollobas}
B{\'e}la Bollob{\'a}s, \emph{Modern graph theory}, GTM, vol. 184,
  Springer-Verlag, New York, 1998.

\bibitem[Cap]{prag}
L.~Caporaso, \emph{Introduction to moduli of curves}, Notes of the Summer
  School PRAGMATIC 2004.

\bibitem[Cap94]{cap}
\bysame, \emph{A compactification of the universal {P}icard variety over the
  moduli space of stable curves}, Journal of the American Mathematical Society
  \textbf{7} (1994), no.~3, 589--660.

\bibitem[Cap08]{capneron}
\bysame, \emph{N\'eron models and compactified {P}icard schemes over the moduli
  stack of stable curves}, Amer. J. Math. \textbf{130} (2008), no.~1, 1--47.

\bibitem[CC03]{CC}
L.~Caporaso and C.~Casagrande, \emph{Combinatorial properties of stable spin
  curves}, Comm. Algebra \textbf{31} (2003), no.~8, 3653--3672, Special issue
  in honor of Steven L. Kleiman.

\bibitem[CCC07]{CCC}
L.~Caporaso, C.~Casagrande, and M.~Cornalba, \emph{Moduli of roots of line
  bundles on curves}, Trans. Amer. Math. Soc. \textbf{359} (2007), no.~8,
  3733--3768 (electronic).

\bibitem[CE02]{C-E}
F.~Chung and R.~B. Ellis, \emph{A chip-firing game and {D}irichlet
  eigenvalues}, Discrete Math. \textbf{257} (2002), no.~2-3, 341--355, Kleitman
  and combinatorics: a celebration (Cambridge, MA, 1999).

\bibitem[Cha82]{chaiken}
S.~Chaiken, \emph{A combinatoric proof of all minors matrix tree theorem}, SIAM
  Journal Algebraic Discrete Methods \textbf{3} (1982), 319.

\bibitem[Chi]{chiodo}
A.~Chiodo, \emph{Quantitative n\'eron theory for torsion bundles},
  arXiv:math.AG/0603689v2.

\bibitem[CR02]{C-R}
H.~Christianson and V.~Reiner, \emph{the {C}ritical {G}roup of a threshold
  graph}, Linear Algebra and its Applications \textbf{349} (2002), 233--244.

\bibitem[CY99]{C-Y}
F.~Chung and S-T. Yau, \emph{Coverings, heat kernels and spanning trees},
  Electron. J. Combin. \textbf{6} (1999), Research Paper 12, 21 pp.\
  (electronic).

\bibitem[DM69]{DM}
P.~Deligne and D.~Mumford, \emph{The irreducibility of the space of curves of
  given genus}, Inst. Hautes \'Etudes Sci. Publ. Math. (1969), no.~36, 75--109.

\bibitem[Gie82]{gie}
D.~Gieseker, \emph{Lectures on moduli of curves}, Tata Institute of Fundamental
  Research Lectures on Mathematics and Physics, vol.~69, Published for the Tata
  Institute of Fundamental Research, Bombay, 1982.

\bibitem[GKM02]{GKM}
A.~Gibney, S.~Keel, and I.~Morrison, \emph{Towards the ample cone of
  {$\overline M\sb {g,n}$}}, J. Amer. Math. Soc. \textbf{15} (2002), no.~2,
  273--294 (electronic).

\bibitem[GM80]{G-McK}
C.~D. Godsil and B.~D. McKay, \emph{Feasibility conditions for the existence of
  walk-regular graphs}, Linear Algebra Appl. \textbf{30} (1980), 51--61.

\bibitem[GR01]{godsil-royle}
C.~D. Godsil and G.~Royle, \emph{Algebraic graphy theory}, Graduate Text in
  Mathematics, no. 207, Springer-Verlag, 2001.

\bibitem[Har69]{Harary}
F.~Harary, \emph{Graph theory}, Addison-Wesley Publishing Co., Reading,
  Mass.-Menlo Park, Calif.-London, 1969.

\bibitem[Har77]{har}
R.~Hartshorne, \emph{Algebraic {G}eometry}, G.T.M., no.~52, Springer-Verlag,
  1977.

\bibitem[HM98]{harris}
J.~Harris and I.~Morrison, \emph{Moduli of {C}urves}, Springer-Verlag, 1998.

\bibitem[JR]{DJ-IR}
D.~Jakobson and I.~Rivin, \emph{On some extremal problems in graph theory},
  preprint arXiv:math/9907050v1 [math.CO].

\bibitem[Kel76]{Kelmans}
A.~K. Kelmans, \emph{Comparison of graphs by their number of spanning trees},
  Discrete Math. \textbf{16} (1976), no.~3, 241--261.

\bibitem[Las01]{bodo}
B.~Lass, \emph{D\'emonstration combinatoire de la formule de {H}arer-{Z}agier},
  C. R. Acad. Sci. Paris S\'er. I Math. \textbf{333} (2001), no.~3, 155--160.

\bibitem[Lor89]{lorarithmetical}
D.~Lorenzini, \emph{Arithmetical graphs}, Math. Ann. \textbf{285} (1989),
  no.~3, 481--501.

\bibitem[Lor90a]{lorgraphs}
\bysame, \emph{Dual graphs of degenerating curves}, Math. Ann. \textbf{287}
  (1990), 135--150.

\bibitem[Lor90b]{lorjac}
\bysame, \emph{Groups of components of {N}\'eron models of {J}acobians},
  Compositio Matematica \textbf{73} (1990), 145--160.

\bibitem[Lor91]{lorfinite}
\bysame, \emph{A finite group attached to the {L}aplacian of a graph}, Discrete
  Math. \textbf{91} (1991), no.~3, 277--282.

\bibitem[Lor93]{lorner}
\bysame, \emph{On the group of components of a {N}\'eron model}, J. Reine
  Angew. Math. \textbf{445} (1993), 109--160.

\bibitem[Lyo05]{lyons}
R.~Lyons, \emph{Asymptotic enumeration of spanning trees}, Combin. Probab.
  Comput. \textbf{14} (2005), 491--522.

\bibitem[McK81]{McKrandom}
B.~D. McKay, \emph{Spanning trees in random regular graphs}, Proceedings of the
  Third Caribbean Conference on Combinatorics and Computing (Bridgetown, 1981)
  (Cave Hill Campus, Barbados), Univ. West Indies, 1981, pp.~139--143.

\bibitem[McK83]{McK}
\bysame, \emph{Spanning trees in regular graphs}, European J. Combin.
  \textbf{4} (1983), no.~2, 149--160.

\bibitem[OS79]{OS}
T.~Oda and C.~Seshadri, \emph{Compactifications of the generalized {J}acobian
  variety}, Trans. A.M.S. \textbf{253} (1979), 1--90.

\bibitem[Ray70]{ray}
M.~Raynaud, \emph{Specialisation du foncteur de {P}icard}, Inst. Hautes Etudes
  Sci. Publ. Math. \textbf{28} (1970), 27--76.

\bibitem[Tsu96]{tsukui}
Y.~Tsukui, \emph{Transformations of {C}ubic {G}raphs}, J. Franklin Institute
  \textbf{333(B)} (1996), no.~4, 565--575.

\bibitem[Tyu03]{tyurin}
A.~Tyurin, \emph{Quantization, classical and quantum field theory and theta
  functions}, CRM Monograph Series, vol.~21, American Mathematical Society,
  Providence, RI, 2003, With a foreword by Alexei Kokotov.

\bibitem[Wes96]{west}
D.~West, \emph{Introduction to graph theory}, Prentice Hall, 1996.

\end{thebibliography}
  
\bigskip

\noindent Margarida Melo, 
Dipartimento di Matematica, Universit\`a di Roma Tre,
Largo S. Leonardo Murialdo 1, 00146 Roma - ITALY\\
E-mail:  \textsl {melo@mat.uniroma3.it}

\bigskip

\noindent Simone Busonero,
Dipartimento di Matematica ``Guido Castelnuovo'', 
Universit\`a di Roma La Sapienza, 
P.le Aldo Moro 2, 00185 Roma - ITALY\\ 
E-mail:  \textsl {busonero@mat.uniroma1.it}

\bigskip
\noindent Lidia Stoppino,
Dipartimento di Matematica, Universit\`a di Pavia, 
Via Ferrata 1, 27100 Pavia - ITALY. \\
E-mail: \textsl {lidia.stoppino@unipv.it}.

\end{document}